\documentclass{article}
\usepackage {amssymb}
\usepackage {amsmath}
\usepackage {amsthm}
\usepackage{graphicx}
\usepackage{subfigure}
\usepackage {amscd}
\usepackage{xypic}
\usepackage[colorlinks, citecolor=red]{hyperref}
\newcommand{\tZ}{{\tilde{Z}}}

\newcommand{\cI}{{\sqrt{-1}}}
\newcommand{\bP}{{\mathbb{P}}}
\newcommand{{\sddbar}}{{\sqrt{-1}\partial\bar{\partial}}}
\newcommand{{\tmfe}}{{\tilde{\mathfrak{e}}}}
\newcommand{{\mfe}}{{\mathfrak{e}}}

\newtheorem{thm}{Theorem}[section]
\newtheorem{prop}{Proposition}[section]

\newtheorem{cor}{Corollary}[section]
\newtheorem{rem}{Remark}[section]

\newtheorem{exmp}{Example}[section]
\newtheorem{lem}{Lemma}[section]

\begin{document}

\title{Numerical Solutions of K\"{a}hler-Einstein metrics on
$\mathbb{P}^2$ with conical singularities along a conic curve}
%\titlerunning{Quantum Monodromy in Integrable Systems}
\author{Chi Li }

%\institute{Mathematics Department, SUNY at Stony Brook, Stony Brook, 11794-3651, USA. \\ \email{chil@math.princeton.edu}}

%\authorrunning{Chi Li}

%\date{Received: 21 April 1998 / Accepted: 8 December 1998}
%\communicated{H. Araki}

\maketitle

\begin{abstract}
We solve for the $SO(3)$-invariant K\"{a}hler-Einstein metric on
$\mathbb{P}^2$ with cone singularities along a smooth conic curve
using numerical approach. The numerical results show the sharp range
of angles ( $(\pi/2,2\pi]$) for the solvability of equations, and
the right limit metric space ($\mathbb{P}(1,1,4)$). These result
exactly match our theoretical conclusion. We also point out the the cause of 
incomplete classifications in \cite{Conti}.
\end{abstract}

\section{Introdution}
Let $D$ be a smooth conic curve in $\mathbb{P}^2$. In this work, we
fix $D=\{Z_1^2+Z_2^2+Z_3^2=0\}$. In the recent work \cite{LiSu}, we
have considered the problem of existence of K\"{a}hler-Einstein
metrics on $\mathbb{P}^2$ with cone singularities along $D$ of cone
angle $2\pi\beta\in (0,2\pi]$. The following is the main result in
this study \cite{LiSu}:
\begin{thm}[\cite{LiSu}]\label{LSmain}
There exists a conical K\"{a}hler-Einstein metric on
$(\mathbb{P}^2,(1-\beta)D)$ if and only if $\beta\in (1/4,1]$.
\end{thm}
As pointed out to us by Dr. H-J. Hein, when $\beta=\frac{1}{3}$,
this gives rise to Calabi-Yau cone metric on the 3-dimensional
$A_2$ singularity $x_1^2+x_2^2+x_3^2+x_4^3=0$.

This is a question raised by Gauntlett-Martelli-Sparks-Yau in
\cite{GMSY}. In \cite{GMSY}, they proved there can not exist such
Calabi-Yau cone metric on 3-dimensional $A_{k-1}$ singularities
$x_1^2+x_2^2+x_3^2+x_4^{k}=0$ if $k\ge 4$. The idea is to look at
the links $L_k$ of such singularities. Any such Calabi-Yau cone
metric would induce a Sasaki-Einstein structure on $L_k$. By further
taking quotient by the $U(1)$ action generated by the natural Reeb
vector field, we would get an orbifold K\"{a}hler-Einstein metric on
$(\mathbb{P}^2, (1-\frac{1}{k})D)$. In \cite{GMSY}, the obstruction
for $k\ge 4$ comes from the Lichnerowics obstruction. In \cite{LiSu}
this was explained as $(\mathbb{P}^2,(1-\frac{1}{k})D)$ being not
log-K-stable if $k\ge 4$. For $k=1$ and $k=2$ case, we have the
standard examples corresponding to the $\mathbb{P}^2$ with
Fubini-Study metric and $(\mathbb{P}^2,\frac{1}{2}D)\cong
\mathbb{P}^1\times\mathbb{P}^1$ with the product metric. These
discussion leaves open the existence problem when $k=3$.

The new insight from \cite{LiSu} is that we can put such kind of
orbifold K\"{a}hler metrics in the more broad family of conical
K\"{a}hler metrics. In our notation $\beta=1/k$. This allows us to
give a uniform theory which together with an interpolation argument
lead us to Theorem \ref{LSmain}.

However, as pointed out in \cite{LiSu}, such result is in
contradiction to the result by Conti in \cite{Conti}, which says
there is no cone Calabi-Yau cone metric on $A_2$ singularities. His
proof is by classifying all the cohomogeneity one 5-dimensional
Sasaki-Einstein manifolds. This leaves us wondering which one is
right.

We decide to attack this question by returning to the approach in
\cite{GMSY} where the equations of orbifold K\"{a}hler-Einstein
metrics on $(\mathbb{P}^2,(1-1/k)D)$ were written down. Note that
because of $SO(3)$ symmetry, such equation comes from the work in
\cite{DaSt}. Moreover, the transformation and change of variables
introduced in \cite{GMSY} is very useful for dealing with the
problem at hand. In this way, we get a 2nd order differential
equation with appropriate boundary conditions.

Since we could not integrate the equation for general $\beta$ we
will use numerical simulation to solve it. This was suggested in
\cite{GMSY}. Our goal is to carry out such numerical approach. As it
turns out, the result is same as we expected.
\begin{thm}\label{main1}
The equations corresponding to $SO(3)$-invariant K\"{a}hler-Einstein
metric $\omega_\beta$ on $(\mathbb{P}^2,(1-\beta)D)$ has a numerical
solution if and only if $\beta>1/4$.
\end{thm}
As suggested by Dr. Song Sun and Dr. H-J. Hein, we will further verify the conjecture
proposed in \cite{LiSu} which predicts the limit metric space as
$\beta$ goes to the critical value $1/4$. Again, the numerical
result fits well with our expectation.
\begin{thm}\label{main2}
As $\beta\rightarrow 1/4$, the metric space
$(\mathbb{P}^2,\omega_\beta)$ converges to the metric space
$(\mathbb{P}(1,1,4),\hat{\omega}_{KE})$ where $\hat{\omega}_{KE}$ is
the induced orbifold K\"{a}hler-Einstein metric coming from the
standard Fubini-Study metric on $\mathbb{P}^2$ by the natural branch
cover: $\mathbb{P}(1,1,1)\rightarrow\mathbb{P}^2(1,1,4)$. Moreover, the bubble out of
this convergence is the $\mathbb{Z}_2$-quotient of Eguchi-Hanson metric on $\mathbb{P}^2\backslash D$. 
\end{thm}
The precise meaning of the above statement is detailed in Section
\ref{numerical} and Section \ref{limit}. These results confirm our result in Theorem \ref{LSmain}. In the last section, we return to calculate
the data of Sasaki-Einstein 5-manifolds associated with $\bP^1\times\bP^1$ and $\bP^2$ in the sense of that in \cite{Conti}. We find that
there are indeed cases ignored in \cite{Conti}.

The example of the pair $(\bP^2,D)$ here can be generalized in more broad settings, which we plan to discuss else where together with Song Sun and H-J. Hein.  

The organization of this note is as follows. The first section gives
a detailed review of the structure of $SO(3)$-orbits for
$\mathbb{P}^2$. The second section discusses the equations we want
to solve. Again, we carefully review the approach in \cite{GMSY} and
work out more details.  In the third sections, we show our first
numerical result Theorem \ref{main1}. In section \ref{limit}, after
describing the $SU(2)$-orbits of $\mathbb{P}(1,1,4)$ we demonstrate
our numerical studies which explains Theorem \ref{main2}. In the last section, we calculate 
the data for $\bP^1\times\bP^1$ in detail. We also calculate the data for the associated Sasaki-Einstein metric which 
indicates the missing case in \cite{Conti}.

\textbf{Acknowledgement}: The author would like to thank Dr. Song Sun and Dr. H-J. Hein for insightful suggestions, and
Professor LeBrun for his interests in our work. The author would also like to thank Dr. Caner Koca for carefully reading the previous version of the paper and pointing out several typos.

\section{$SO(3)$ orbits}\label{orbitsP2}
Let us first review how to decompose $\mathbb{P}^2$ into
$SO(3,\mathbb{R})$-orbits following \cite{GMSY}. First note that
$\mathbb{P}^2=(\mathbb{C}^3-\{0\})/\mathbb{C}^*$ under the
equivalence relation $(Z_1,Z_2,Z_3)\sim (\lambda Z_1,\lambda
Z_2,\lambda Z_3)$ for some $\lambda\neq 0\in\mathbb{C}^*$. Now fix
any $0\neq Z:=(Z_i)_{i=1}^3\in \mathbb{C}^3$, it determines a point
in $\mathbb{P}^2$ with homogeneous coordinate
$[Z]:=[Z_i]_{i=1}^3=[Z_1,Z_2,Z_3]$. Now write the polar
decomposition
\[
Z_1^2+Z_2^2+Z_3^2=\rho^2 e^{2i\theta}.
\]
So if we define
\[
\tZ_i= e^{-i\theta}Z_i,
\]
then $[Z_i]_{i=1}^3=[\tZ_i]_{i=1}^3$ and
\begin{equation}\label{normalize}
\tZ_1^2+\tZ_2^2+\tZ_3^2=\rho^2\ge 0.
\end{equation}
Now write
\[
\tZ_i=u_i+\sqrt{-1}v_i,
\]
then the identity \eqref{normalize} is
equivalent to the identity
\begin{equation}
|u|^2-|v|^2=\rho^2;\quad u\cdot v=0.
\end{equation}
We use these two relations to define the set:
\[
\mathbb{O}=\{(u,v):=u+iv\neq 0|u\cdot v=0, |u|^2-|v|^2\ge 0
\}\subset(\mathbb{R}^3)^2-\{0\}.
\]
Define an equivalence relation on $\mathbb{O}$ by \footnote{Dr. Caner Koca pointed out to me that in the second case, the multiplication of $e^{i\theta}$ was missing in the previous version of the paper.}
\[
\left\{\begin{array}{lll} (u,v)\sim a(u,v)&, \forall a\in\mathbb{R}^{\times}&, \mbox{ if } |u|\neq |v|;\\
(u,v)\sim a e^{i\theta}(u,v)&, \forall a\in \mathbb{R}^{\times}, \forall \theta\in [0,2\pi)&,  \mbox{ if }  |u|=|v|. \end{array}\right.
\]
%$(u,v)\sim a(u,v)$
%for some $a\in\mathbb{R}^{\times}$. 
Denote the quotient set by
$\overline{\mathbb{O}}=\mathbb{O}/\sim$. Then we have defined a
homeomorphism
\begin{eqnarray*}
\Phi: \mathbb{P}^2&\longrightarrow& \overline{\mathbb{O}}\\
{[Z_i]_{i=1}^3}&\mapsto& [u,v] \mbox{ satisfying
}u+\sqrt{-1}v=e^{-\frac{i}{2}
\rm{Arg}(Z_1^2+Z_2^2+Z_3^2)}(Z_1,Z_2,Z_3).
\end{eqnarray*}
Here we assume ${\rm Arg}(0)$ can be any real number, which is compatible with the 2nd case in the equivalence.
The $SO(3)$ acts on $\mathbb{P}^2\cong\overline{\mathbb{O}}$ by
\[
g\cdot (u,v)=(gu,gv).
\]
The quotient of this action is an interval:
\begin{eqnarray*}
R:\overline{\mathbb{O}}&\longrightarrow& [0,1]\\
{[u,v]}&\mapsto& \frac{|v|}{|u|}
\end{eqnarray*}
So the function $R$ classifies $SO(3)$ orbit. Moreover it's easy to
verify that equivalently we have the relation
\begin{equation}\label{Rinv}
\frac{|Z_1^2+Z_2^2+Z_3^2|}{|Z_1|^2+|Z_2|^2+|Z_3|^2}=\frac{1-R^2}{1+R^2}.
\end{equation}

For each point $(u,v)\in\mathbb{O}$, we get an orthonormal basis in the following way. If $v\neq 0$, we set $(e_u=u/|u|,e_v=v/|v|,e_w:=e_u\times e_v)$. If $v=0$ We choose any $e_v$ perpendicular to $e_u=u/|u|$ and let $e_w=e_u\times u_v$. We will denote
$U(1)_1$, $U(1)_2$ and $U(1)_3$ to be the rotation around the axes
in the direction $e_u$, $e_v$ and $e_w$ respectively.
\begin{lem}\label{orbits}
The generic orbit is ${\rm Orb}_{R=R_0}=SO(3)/\mathbb{Z}_2$ (when
$0<R_0=R([u,v])<1$). The two special orbits are
\begin{eqnarray*}
{\rm Orb}_{R=0}&=&(SO(3)/\mathbb{Z}_2)/U(1)_1=\mathbb{RP}^2; \\
{\rm Orb}_{R=1}&=&(SO(3)/\mathbb{Z}_2)/U(1)_3=\mathbb{P}^1.
\end{eqnarray*}
\end{lem}
\begin{proof}
When $0<R=\frac{|v|}{|u|}<1$, the stabilizer of $SO(3)$ action at
$[v,w]$ is isomorphic to $\mathbb{Z}_2$ with generator being the
rotation around $e_w$ with angle $\pi$, i.e.
$(e_u,e_v,e_w)\rightarrow (-e_u,-e_v,e_w)$.

When $R=0$, $v$=0. The stabilizer is generated by $\mathbb{Z}_2$ and
$U(1)_1$. The generator of $\mathbb{Z}_2$ can be chosen to be
$(e_u,e_v,e_w)\mapsto (-e_u,-e_v,e_w)$ (for any $e_v$, $e_w$ such
that $\{e_u,e_v,e_w\}$ is an orthonormal basis). $U(1)_1$ is the
rotation group around $e_u$. It's easy to verify that
\[
{\rm Orb}_{R=0}=(\mathbb{R}^3-\{0\})/\mathbb{R}^{\times}=\mathbb{RP}^2.
\]
When $R=1$, $|u|=|v|$. The stabilizer is $U(1)$-rotation group
around $e_w$ denoted as $U(1)_3$. Note $\mathbb{Z}_2\subset U(1)_3$.
It's easy to see that (for example by \eqref{Rinv})
\[
{\rm
Orb}_{R=1}=\{Z_1^2+Z_2^2+Z_3^2=0\}\cong\mathbb{P}^1\subset\mathbb{P}^2.
\]
\end{proof}
Fix the generator of $so(3)=\rm{Lie} (SO(3))$ to be
\[
X_1=\left(\begin{array}{ccc} 0&0&0\\0&0&-1\\0&1&0
\end{array}\right)\;, X_2=\left(\begin{array}{ccc} 0&0&1\\0&0&0\\-1&0&0
\end{array}\right),\; X_3=\left(\begin{array}{ccc} 0&-1&0\\1&0&0\\0&0&0
\end{array}\right).
\]
Then the corresponding invariant vector field on the orbit
$SO(3)([u,v])$ at point $[u,v]$ is given by the infinitesimal
rotation around three axes in the directions of $e_u$, $e_v$, $e_w$
respectively. In other words, they are generators of the action
of $U(1)_1$, $U(1)_2$, $U(1)_3$ respectively.
\begin{enumerate}
\item Around $e_u$:
\begin{eqnarray*}
T_u&=&\left.\frac{d}{d\theta}\right|_{\theta=0}(u+\sqrt{-1}(\cos\theta
e_v-\sin\theta e_w)|v|)=-\sqrt{-1} |v|e_w.
\end{eqnarray*}

\item Around $e_v$:
$ T_v=\left.\frac{d}{d\theta}\right|_{\theta=0}(\sin\theta
e_w+\cos\theta e_u)|u|+\sqrt{-1}v=|u|e_w$.

\item Around $e_w$:
\begin{eqnarray*}
T_w&=&\left.\frac{d}{d\theta}\right|_{\theta=0}(|u|(\cos\theta
e_u-\sin\theta
e_v)+\sqrt{-1}(\sin\theta e_u+\cos\theta e_v)|v| )\\
&=& -|u| e_v+\sqrt{-1} |v|e_u.
\end{eqnarray*}
\end{enumerate}
We can define another vector field generating the radial
transformation
\[
T_R=\left.\frac{d}{d\theta}\right|_{\theta=0}\left(|u|(e_u+\sqrt{-1}\left(\frac{|v|}{|u|}+\theta\right)e_v\right)=\sqrt{-1}|u|e_v.
\]
Note that the above vectors represent the tangent vector in
\[T_{[u+iv]}\mathbb{P}^2={\rm
Hom}(\mathbb{C}(u+iv),(\mathbb{C}(u+iv))^{\perp})\cong {\rm
Hom}(\mathbb{C}(u+iv),\mathbb{C}^{3}/\mathbb{C}(u+iv)).\] 
\begin{lem}\label{vanish}
On ${\rm Orb}_{R=0}=\mathbb{RP}^2$, $T_u=0$; On ${\rm
Orb}_{R=1}=\mathbb{P}^1$, $T_w=0$.
\end{lem}
\begin{proof}
When $R=0$, $|v|=0$, so $T_u=0$ on ${\rm Orb}_{R=0}RP^2$. When
$R=1$,
\[
T_w=|u|\frac{v}{|v|}-\sqrt{-1}|v|\frac{u}{|u|}=-\sqrt{-1}(u+\sqrt{-1}v)
\]
so $T_w=-\sqrt{-1}(u+\sqrt{-1}v)\in\mathbb{C}\cdot(u,v)$, so
$T_w|_{R=1}=0$, i.e. $T_w$ vanishes on the special orbits ${\rm
Orb}_{R=1}=\mathbb{P}^1$. .
\end{proof}
Note that this Lemma also follows from Lemma \ref{orbits} by the
fact that $U(1)_1$ is the stabilizer group on ${\rm Orb}_{R=0}$
generated by $T_u$, while $U(1)_3$ is the stabilizer group on ${\rm
Orb}_{R=1}$ generated by $T_w$.

\section{Equations for $SO(3)$ invariant K\"{a}hler-Einstein}\label{secKEeq}

For special metrics $g$ on $\mathbb{P}^2$, we have the following
\begin{lem}\label{restrictabc}
\begin{enumerate}
\item
For any K\"{a}hler metric $g$, we have $|T_u|_g\le |T_v|_g$. The
equality holds only on the special orbit ${\rm
Orb}_{R=1}=\mathbb{P}^1$.
\item
For any $SO(3)$ invariant metric $g$, $|T_v|_g=|T_w|_g$ on the
special orbit ${\rm Orb}_{R=0}=\mathbb{RP}^2$.
\end{enumerate}
\end{lem}
\begin{proof}
\begin{enumerate}
\item
Because K\"{a}hler metric is compatible with complex structure
$J=i\cdot$, so
\begin{equation}\label{aleb}
0\le\frac{|T_u|_g}{|T_v|_g}=\frac{|i|v|e_w|_g}{||u|e_w|_g}=\frac{||v|e_w|_g}{||u|e_w|_g}=\frac{|v|}{|u|}=R\le
1.
\end{equation}
\item On the special orbit ${\rm Orb}_{R=0}=\mathbb{RP}^2$, $v=0$.
Let $\gamma_1(\theta)=|u|(-\sin\theta e_w+\cos\theta e_u)$ and
$\gamma_2(\theta)=|u|(\cos\theta e_u+\sin\theta e_v)$. Then
$T_v=\gamma_1'(0)$ and $T_w=\gamma_2'(0)$. Because there exist
rotations $g(\theta)$ in $SO(3)$ such that $g(\theta)\cdot
\gamma_1(\theta)=\gamma_2(\theta)$, the conclusion follows from
invariance of the metric under $SO(3)$.
\end{enumerate}
\end{proof}
Now choose the dual basis of $\{T_R,T_u,T_v,T_w\}$ to be one forms
given by $\{dR,\sigma_1,\sigma_2,\sigma_3\}$. For any $SO(3)$
invariant K\"{a}hler  metric on $\mathbb{P}^2$,
$\{T_u,T_v,T_w,T_R\}$ is orthogonal. The metric can be written in
the form
\begin{equation}\label{so(3)form}
g=(dt)^2+a^2\sigma_1^2+b^2\sigma_2^2+c^2\sigma_3^2.
\end{equation}
where
\[
dt=-|T_R|_gdR,\quad a=|T_u|_g, \quad b=|T_v|_g, \quad c=|T_w|_g.
\]
The minus sign in the first identity is to make the special orbit
$\mathbb{P}^1$ to sit in the distance 0 location. By Lemma
\ref{vanish} and Lemma \ref{restrictabc}, we know that
\begin{cor}\label{corabc} For any $SO(3)$-invariant K\"{a}hler metric on $\mathbb{P}^2$, 
we have $a\le b$ on $\mathbb{P}^2$. On ${\rm Orb}_{R=1}=\mathbb{P}^1$,
$c=0$, $a=b$. On ${\rm Orb}_{R=0}=\mathbb{RP}^2$, $a=0$, $b=c$.
\end{cor}
\begin{exmp}\label{P2FS}
When $\beta=1$, then the $SO(3)$ invariant metric is the
standard Fubini-Study metric on $\mathbb{P}^2$. We can write it in
the form of \eqref{so(3)form}. One way to do this is to recall the
following description of Study-Fubini metric. Let
$\gamma(t):=[Z_1(t),Z_2(t),Z_3(t)]$ be a curve in $\mathbb{P}^2$
with the tangent vector is
$\gamma'(0)=((Z_1(0),Z_2(0),Z_3(0))\mapsto
(Z_1'(0),Z_2'(0),Z_3'(0)))\in {\rm
Hom}(\mathcal{O}_{\mathbb{P}^2}(1),\mathbb{C}^{3}/\mathcal{O}_{\mathbb{P}^2}(1))$.
The length of $\gamma'(0)$ is given by
\[
|\gamma'(0)|_{FS}^2=\frac{|Z'(0)^{\perp}|^2}{|Z(0)|^2}=\left(|Z'(0)|^2-\frac{|\langle
Z'(0),Z(0)\rangle|^2}{|Z(0)^2|}\right)/|Z(0)|^2,
\]
where $\langle\cdot,\cdot\rangle$ is the standard real inner product on $\mathbb{C}^3\cong\mathbb{R}^6$. Using this formula, it's easy to verify that
\begin{eqnarray*}
&|T_R|_{FS}=\frac{|u|^2}{|u|^2+|v|^2}=\frac{1}{1+R^2},\quad
&|T_u|_{FS}=\frac{|v|}{\sqrt{|u|^2+|v|^2}}=\frac{R}{\sqrt{1+R^2}}\\
&|T_v|_{FS}=\frac{|u|^2}{\sqrt{|u|^2+|v|^2}}=\frac{1}{\sqrt{1+R^2}},\quad
&|T_w|_{FS}=\frac{|u|^2-|v|^2}{|u|^2+|v|^2}=\frac{1-R^2}{1+R^2}.
\end{eqnarray*}
So the normal distance function $t$ is determined by
\[
dt=-\frac{1}{1+R^2}dR \quad \&\quad t(1)=0 \Longrightarrow
R=\tan\left(\frac{\pi}{4}-t\right).
\]
So $0\le t\le \pi/4$ and
\[
a=\sin\left(\frac{\pi}{4}-t\right)=\cos\left(t+\frac{\pi}{4}\right),\quad
b=\sin\left(t+\frac{\pi}{4}\right),\quad
c=\cos\left(\frac{\pi}{2}-2t\right)=\sin(2t).
\]
\end{exmp}
\begin{exmp}\label{P1P1FS}
The data for $\bP^1\times\bP^1=(P^2,\frac{1}{2}D)$ are given as
follows. See section \ref{p1p1} for the derivation of these data.
(See also \cite{DaSt} and \cite{GMSY})
\[
a(t)=\frac{1}{\sqrt{3}}\cos(\sqrt{3}t),\quad
b(t)=\frac{1}{\sqrt{3}},\quad c(t)=\frac{1}{\sqrt{3}}\sin
(\sqrt{3}t).
\]
The range for $t$ is $0\le t\le\pi/(2\sqrt{3})$.
\end{exmp}
By \cite{DaSt} and \cite{GMSY}, the equation for K\"{a}hler-Einstein with \textbf{ Ricci curvature equal to 6 } is
reduced to a system of ODEs:
\begin{equation}\label{KEodes}
\left\{
\begin{array}{ccl}
\dot{a}&=&-\frac{b^2+c^2-a^2}{2bc}\\
&&\\
\dot{b}&=&-\frac{a^2+c^2-b^2}{2ac}\\
&&\\
\dot{c}&=&-\frac{a^2+b^2-c^2}{2ab}+6ab
\end{array}\right. 0\le t\le t_*=t_{max}.
\end{equation}
Note that the equation in \cite{GMSY} defers from \cite{DaSt} by a
(negative) factor ($-abc$) which is caused by a change of variable.

The boundary condition at $t=0$ corresponds to the special
orbit ${\rm Orb}_{R=1}=\mathbb{P}^1$ where by Corollary \ref{corabc}
$a=|T_u|_g=|T_v|_g=b$ and $c=|T_w|_g=0$. Moreover, the cone angle
equal to $2\pi\beta$ along ${\rm Orb}_{t=0}=\mathbb{P}^1$ requires
$\dot{c}=2\beta$. The factor $2$ comes from the fact that when
$0<R<1$ the stabilizer is $\mathbb{Z}_2$. So the boundary is given
\begin{eqnarray*}
a(t)&=&\alpha+O(t)\\
b(t)&=&\alpha+O(t)\\
c(t)&=&2\beta t+O(t^2)
\end{eqnarray*}
Since the normalized K\"{a}hler-Einstein metric $\omega'_\beta$
satisfies
\[
Ric(\omega'_\beta)=3\omega'_\beta+2\pi(1-\beta)\{D\}.
\]
Because $[D]=\frac{2}{3} c_1(\mathbb{P}^2)$, so, by taking
cohomological classes on both sides, we get
\[
3[\omega'_\beta]=\frac{1}{3}(1+2\beta)\cdot 2\pi c_1(\mathbb{P}^2).
\]
So $\alpha$ and $\beta$ are related by $ \alpha^2=\delta\cdot
\frac{1}{3}(1+2\beta)$ since both sides are proportional to the
volume of $\mathbb{P}^1$. The factor $\delta$ can be carefully
tracked out, but it can also be easily determined either by checking
the standard $\mathbb{P}^2$ with Fubini-Study metric in Example
\ref{P2FS} or by substituting in to the last equation in
\eqref{KEodes}. The result is
\[
\alpha^2=\frac{1}{6}(1+2\beta).
\]
%Note that the relation between $\alpha$ and $\beta$ can also be
obtained from equation \eqref{KEodes}.

When $t=t_*=t_{max}$, we know from Corollary \ref{corabc} that $a(t_*)=0$
and $b(t_*)=c(t_*)$.
\begin{lem}
\begin{equation}\label{dert*}
\dot{a}(t_*)=-1,\quad \dot{b}(t_*)=\dot{c}(t_*)=0.
\end{equation}
\end{lem}
\begin{proof}
From the first equation in \eqref{KEodes} and $b(t_*)=c(t_*)$, we
get $\dot{a}(t_*)=-1$. Then we use this to derive from Equation
\eqref{KEodes} that
\[
\dot{b}(t_*)=-\dot{c}(t_*)=\lim_{t\rightarrow
t_*}\frac{b-c}{a}=-(\dot{b}(t_*)-\dot{c}(t_*))=-2\dot{b}(t_*).
\]
So the 2nd identity follows.
\end{proof}
Note that $\dot{a}(t_*)=-1$ is compatible with the fact that the
metric is smooth along ${\rm Orb}_{R=0}\cong\mathbb{RP}^2$.

Note the solutions of equation \eqref{KEodes} is not unique around
the point $(a(0),b(0),c(0))=(\alpha,\alpha,0)$. There are at least
three possibilities: $a\le b$, $a=b$, $a\ge b$. The $a=b$ case
corresponds to the Gibbons-Pope-Pederson metric as pointed out in
\cite{DaSt}. We are in the $a\le b$ case. The symmetry of $a$, $b$
is broken by writing down the differential equation for the variable
$R=a/b$. Using \eqref{KEodes}, we get
\[
c\frac{d}{dt}\left(\frac{a}{b}\right)=\left(\frac{a}{b}\right)^2-1.
\]
So it's natural to do the following change of variables introduced
by \cite{GMSY}.
\begin{equation}\label{coordtr}
\frac{dr}{dt}=1/c.
\end{equation}
Then
\[
\frac{dR}{dr}=R^2-1.
\]
Using $a\le b$ (\eqref{aleb}), we get the solution
\begin{equation}\label{aoverb}
R=\frac{a}{b}=-\tanh (r).
\end{equation}
Moreover, we get the range for $r$: $-\infty<r\le 0$. We list the the ranges of $R,t,r$ as follows:\\
\begin{center}
\begin{tabular}{c|c|c|c}
&$\mathbb{P}^2$& $SO(3)/\mathbb{Z}_2$& $RP^2$\\
\hline
$R$&$R=1$&$1>R>0$& $R=0$\\
\hline
$t$&$t=0$&$0<t<t_*$& $t=t_*$\\
\hline $r$&$r=-\infty$&$-\infty<r<0$& $r=0$
\end{tabular}
\end{center}

Define $f=ab$, then $f$ satisfies the second order differential
equation 
\[
\frac{d}{dr}\log\left(f\frac{df}{dr}\right)=2[6f+\coth(2r)].
\]
\begin{exmp}
By easy calculations, one can get that, for $\mathbb{P}^2$,
$f=-\frac{1}{2}\tanh(2r)$,  $f_r(0)=-1$; and for
$\mathbb{P}^1\times\mathbb{P}^1$, $f=-\frac{1}{3}\tanh(r)$,
$f_r(0)=-\frac{1}{3}$. See \cite{GMSY} and also Section \ref{p1p1}.
\end{exmp}
Let $h=f_r$ then this is equivalent to a system:
\begin{equation}\label{fgodes}
\left\{
\begin{array}{ccl}
f_r&=&h\\
&&\\
h_r&=&12 fh+2\coth(2r)h-\frac{h^2}{f}.
\end{array}\right.
\end{equation}
It's easy to verify that the data $(f,R,h)$ and $(a,b,c)$ determine
each other by the relation
\begin{equation}\label{2sets}
f=ab,\quad R=\frac{a}{b},\quad h=f_r=-c^2.
\end{equation}
The boundary condition is given by
\[
f(-\infty)=\alpha^2, f(0)=0.
\]
\[
h(0)=f_r(0)=-c(t_*)^2=-b(t_*)^2.
\]
Using \eqref{dert*}, \eqref{KEodes} and $t_r(t)=c(t)$, we get
\[
h_r(0)=f_{rr}(0)=\left((f_{tt}t_r+f_t(t_r)_t)t_r\right)|_{t=t_*}=\ddot{a}(t_*)b(t_*)^3=0.
\]

\section{Numerical Studies: $\beta>1/4$}\label{numerical}
\begin{figure}[h]
  \begin{center}
    \subfigure[$(r,f)$]{\label{P2-a}\includegraphics[height=4cm]{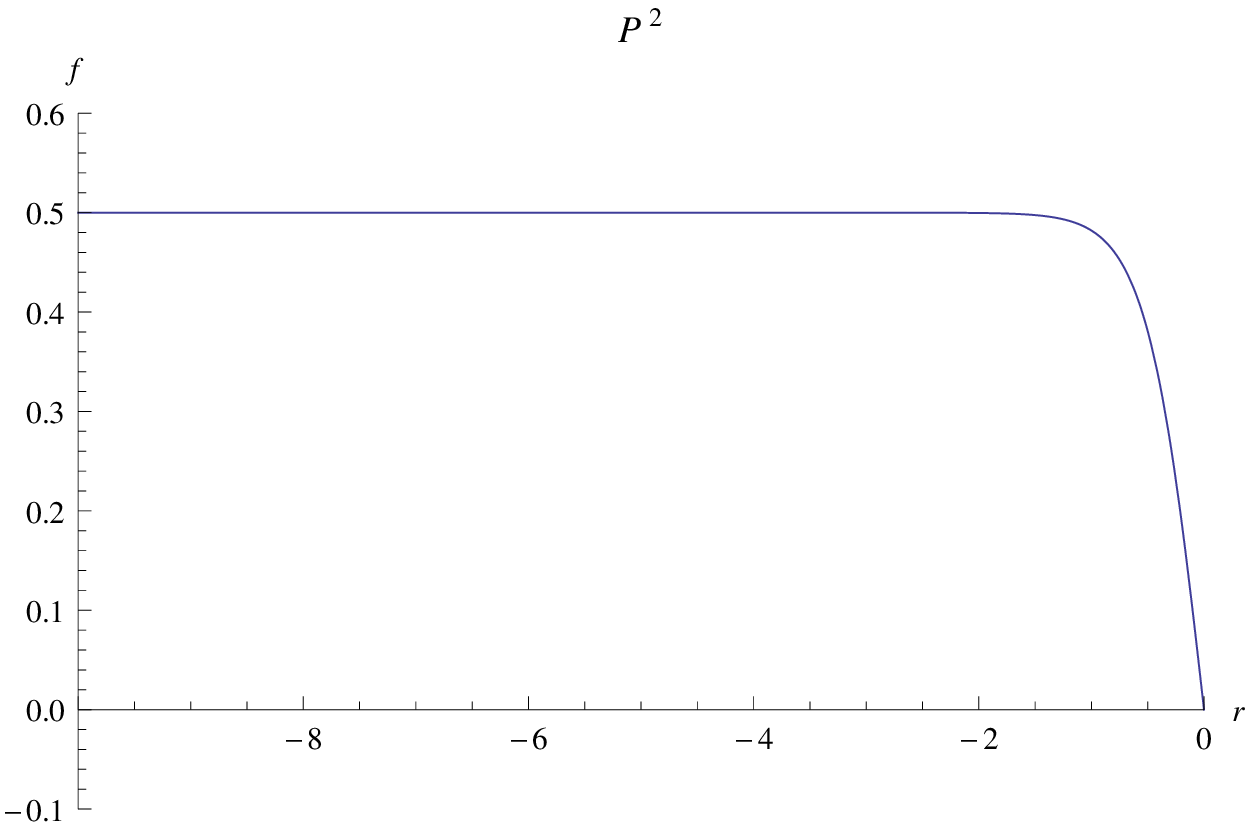}}
    \subfigure[$(r,f_r)$]{\label{P2-b}\includegraphics[height=4cm]{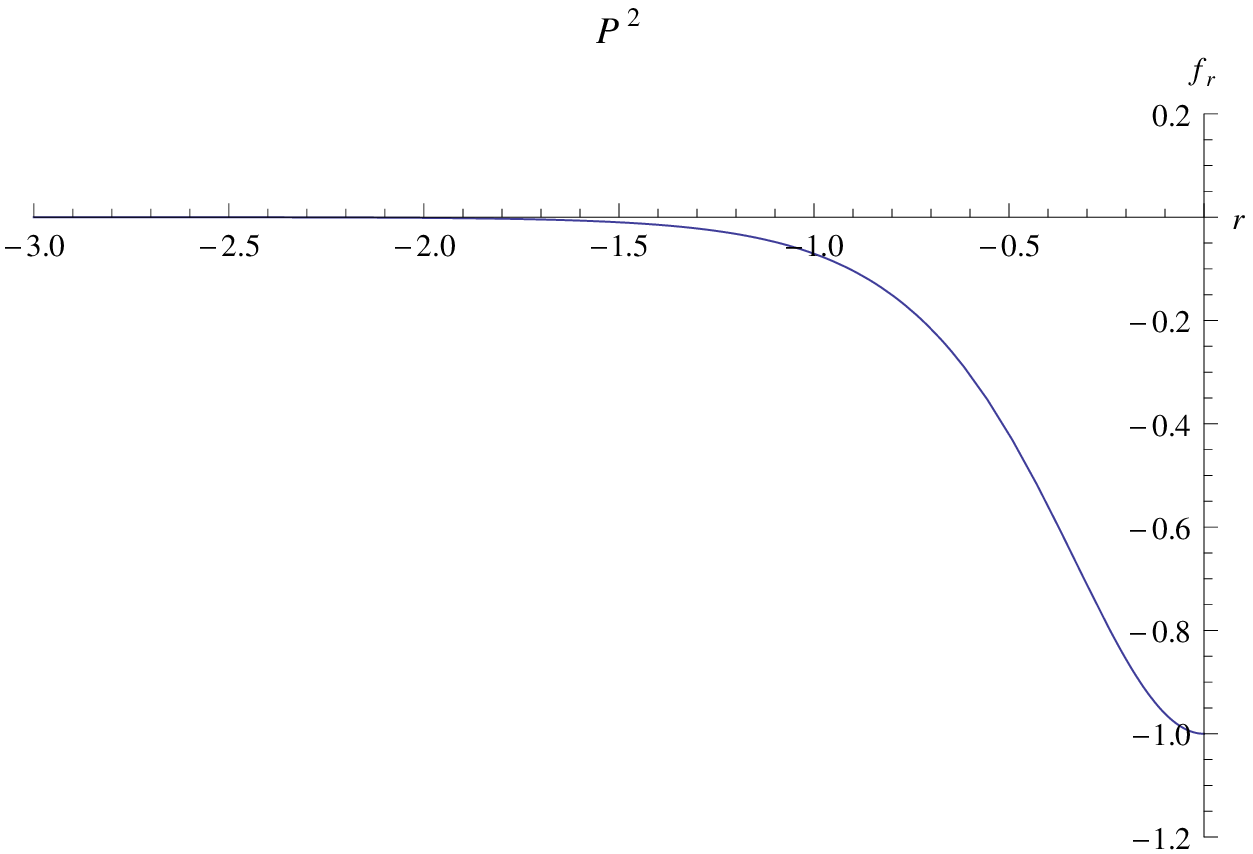}}
  \end{center}
  \caption{Data for $\mathbb{P}^2$}
  \label{P2}
\end{figure}

\begin{figure}[h]
  \begin{center}
    \subfigure[$(r,f)$]{\label{P1P1-a}\includegraphics[height=4cm]{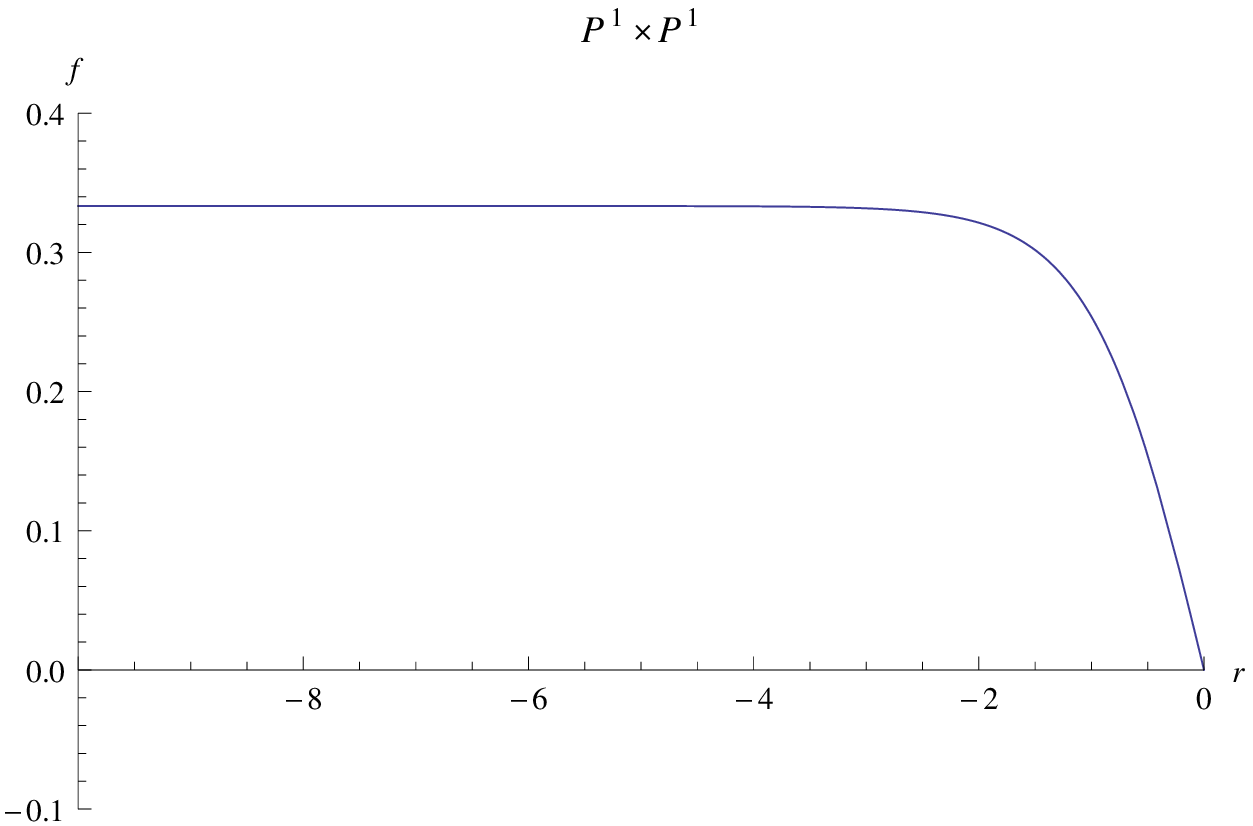}}
    \subfigure[$(r,f_r)$]{\label{P1P1-b}\includegraphics[height=4cm]{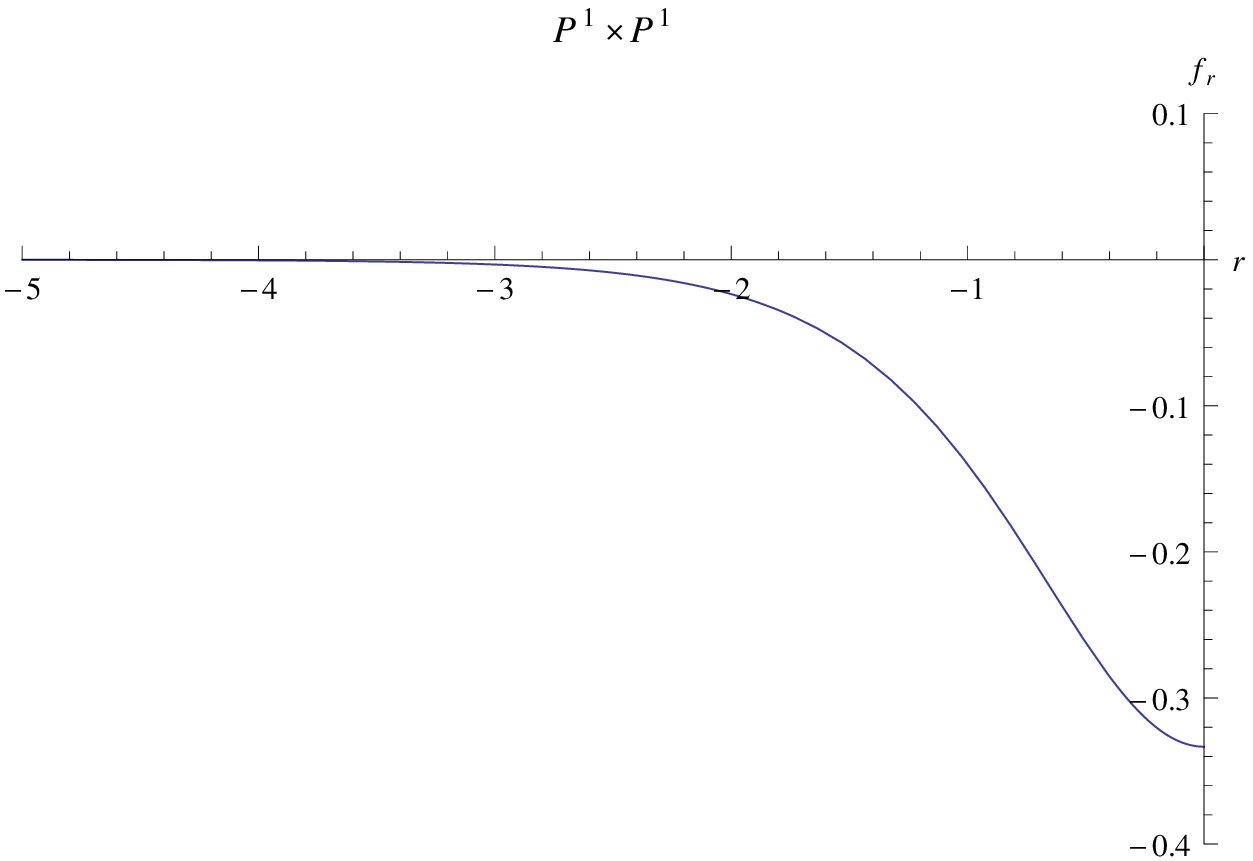}}
  \end{center}
  \caption{Data for $\mathbb{P}^1\times\bP^1$}
  \label{P1P1}
\end{figure}

Now we explain our numerical simulation. We introduce the variable
$\tau$ for convenience and choose boundary value
$(f(0),h(0))=(0,-\frac{1}{\tau}:=-b(t_*)^2)$ and solve the equation
\eqref{fgodes} numerically. However, this can not be done because
there is a zero on the denominator for $r=0$ on the second equation
in \eqref{fgodes} (although it's cancelled by zero on the numerator).
We can however move away from $r=0$ a little bit by using the
boundary condition and Taylor expansion:
\begin{eqnarray*}
f(r)&=&f(0)+f_r(0)r+O(r^2)=-\frac{1}{\tau}r+O(r^2)\\
h(r)&=&h(0)+h_r(0)r+O(r^2)=-\frac{1}{\tau}+O(r^2)
\end{eqnarray*}
So numerically, we can choose $r_0<0$ to be very close to $0$ and
choose the boundary condition to be
\[
(f(r_0),h(r_0))=(-\frac{r_0}{\tau},-\frac{1}{\tau}).
\]
For example, in the following numerical simulation, we choose
$r_0=-10^{-5}$. Then we can shoot the trajectory out for $r$ going
from $r_0$ backward to $-\infty$. Figure \ref{P2} and Figure \ref{P1P1} are the numerical solution
corresponding to $\mathbb{P}^2$ when $\tau=1$ and
$\mathbb{P}^1\times\mathbb{P}^1=(P^2,\frac{1}{2}D)$ when $\tau=3$
respectively. They can be obtained for example by the
\textbf{NDSolve} tool in \textbf{Mathematica}.

%\begin{figure}
%\begin{center}
%\caption{P2}
%\includegraphics[height=4cm]{P2_1.eps}
%\end{center}
%\end{figure}
%\includegraphics[height=4cm]{P2_2.eps}
%\includegraphics[height=4cm]{P1P1_1.eps}
%\includegraphics[height=4cm]{P1P1_2.eps}
%\end{center}
Of course, the above graphs of $f=f(r)$ just recover the graph
$f(r)=-\frac{1}{2}\tanh(2r)$ for $\mathbb{P}^2$ and
$f(r)=-\frac{1}{3}\tanh(r)$ for $\mathbb{P}^1\times\mathbb{P}^1$ (up
to high precision).

If we choose different $\tau$, then we get different solution $f$,
$h=f_r$. We know that $\lim_{r\rightarrow-\infty}
f(r)=f(-\infty)=\alpha^2=\frac{1+2\beta}{6}$. Numerically, we can
just evaluate $f(r)$ for $r$ being sufficiently negative to
calculate $\alpha^2$. Actually, after several tests, one can observe
that for fixed $\tau$, the graph will becomes flat as $r$ goes toward $-\infty$ which
means $f(r)$ becomes stabilized. The speed of approaching flatness depends on the boundary 
value $h(0)=-\frac{1}{\tau}$. The bigger $\tau$ is, the longer $r$-distance it takes for the graph to become flat. 
(This is related to the bubbling phenomenon below)

\begin{figure}[h]
  \begin{center}
    \subfigure[short range]{\label{taualpha-a}\includegraphics[height=4cm]{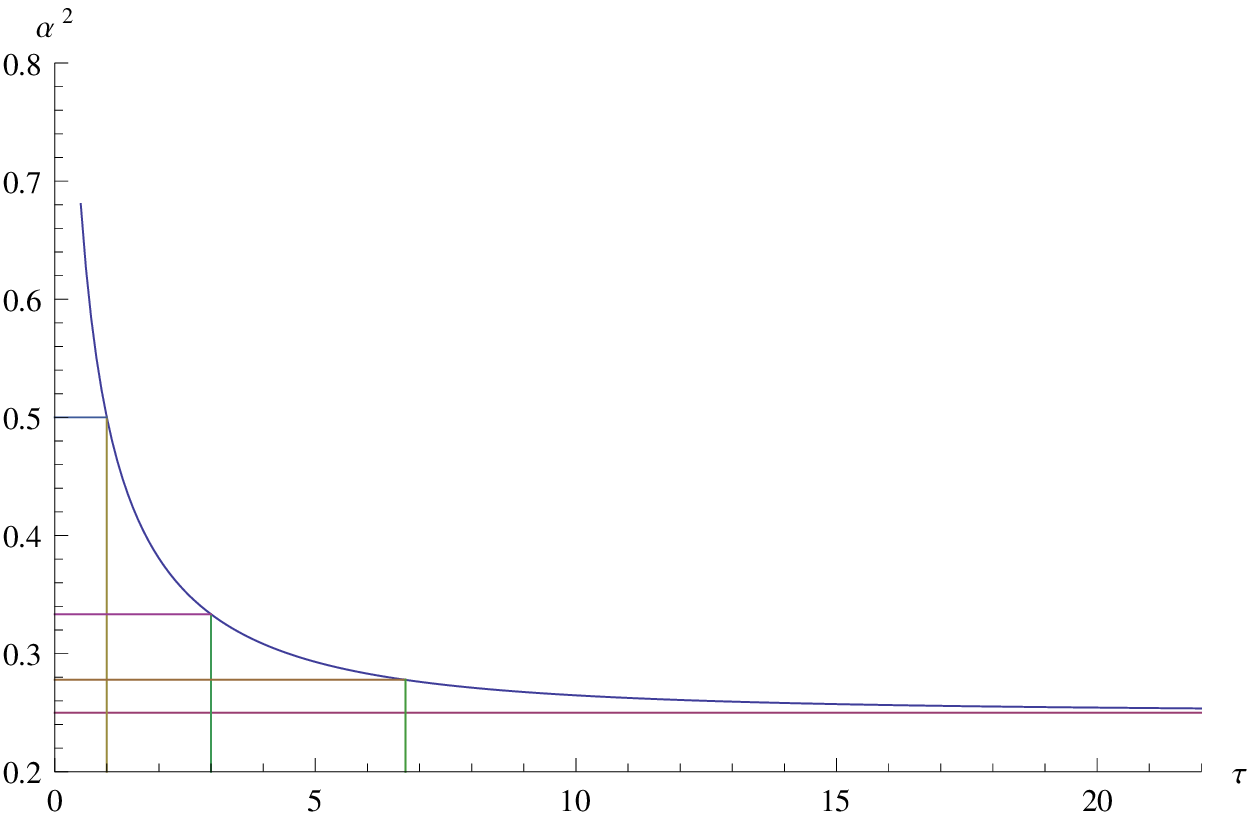}}
    \subfigure[long range]{\label{taualpha-b}\includegraphics[height=4cm]{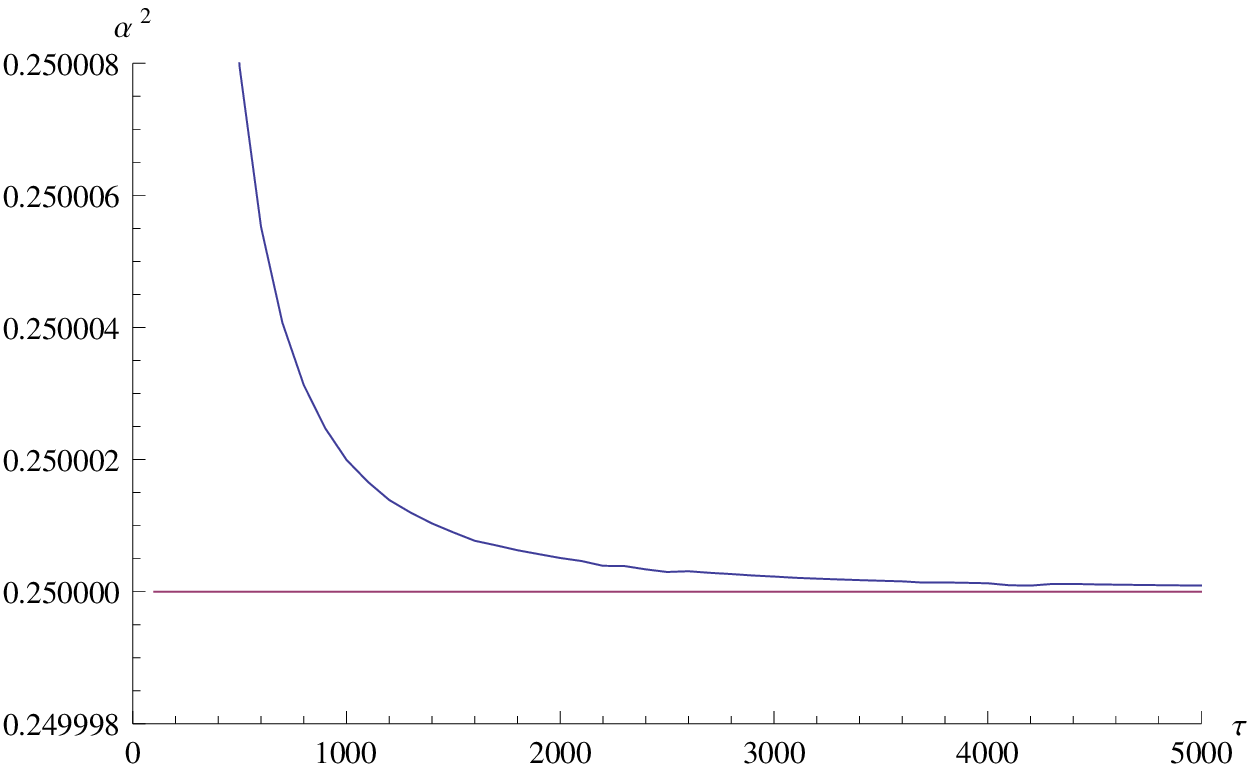}}
  \end{center}
  \caption{$(\tau, \alpha^2)$}
  \label{taualpha}
\end{figure}

We can use \textbf{Mathematica} to calculate (very dense)
sequences of data for $\{\tau, f(\tau,r)\}$ where we make solution
$f$ depend the boundary data $\tau$. Then we sample the value of
$f(\tau,r)$ at $r=-500$. (One can certainly choose $r$ to be more
negative but the visual effect does not change) Figure \ref{taualpha} shows the
numerical result. The two subfigures are for short range and long
range of $\tau$ respectively.
%\begin{center}
%\includegraphics[height=4cm]{data2.eps}
%\includegraphics[height=4cm]{longta.eps}
%\end{center}

\begin{figure}[h]
  \begin{center}
    \subfigure[$(r,f)$]{\label{P23D-a}\includegraphics[height=4cm]{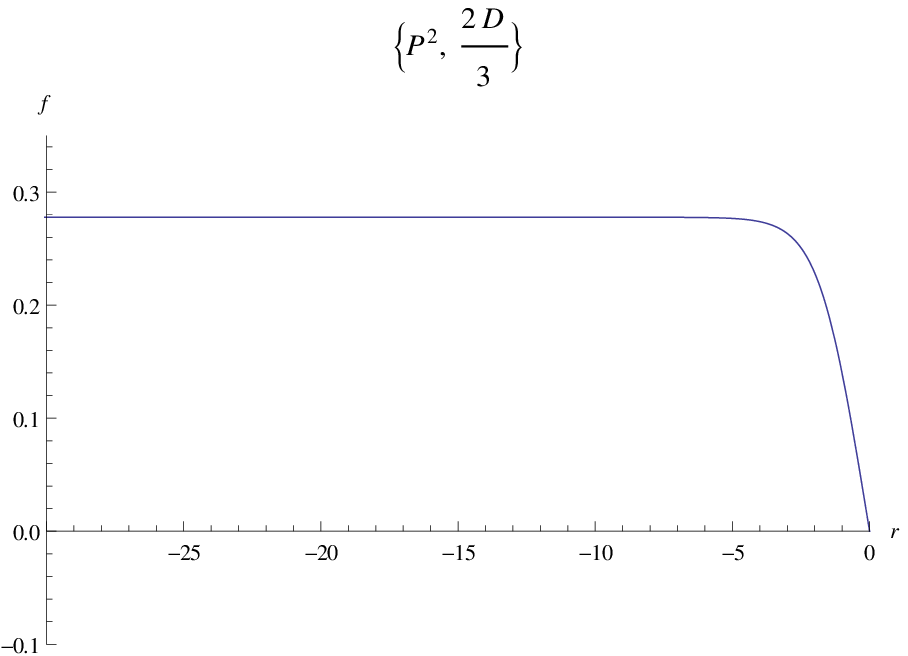}}
    \subfigure[$(r,f_r)$]{\label{P23D-b}\includegraphics[height=4cm]{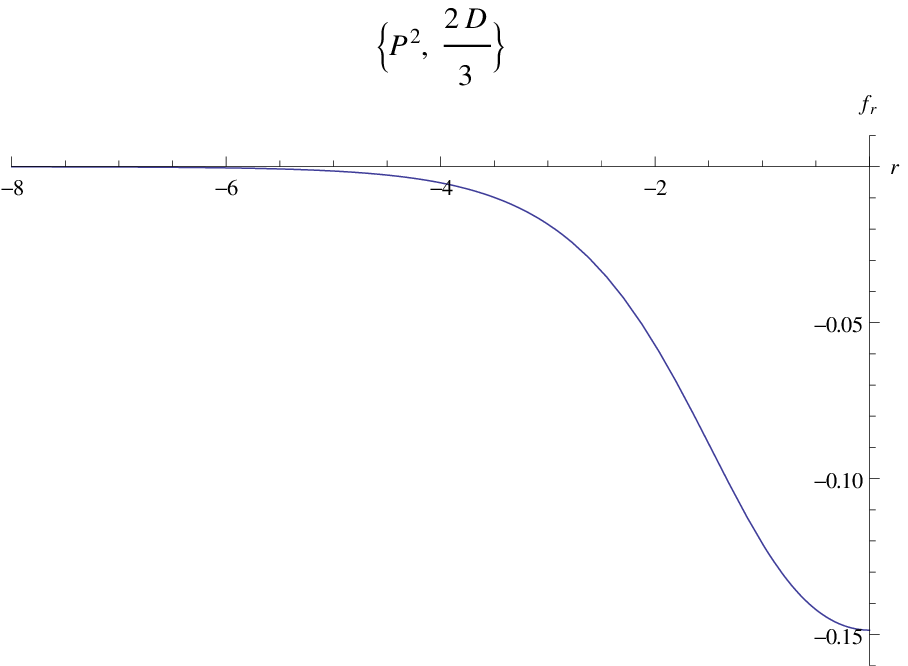}}
  \end{center}
  \caption{Data for $(\bP^2,\frac{2}{3}D)$}
  \label{P23D}
\end{figure}

We see immediately that $\alpha^2$ is a decreasing function of
$\tau$. More importantly, from the picture, we sees that one always
has
\[
\alpha^2=\frac{1+2\beta}{6}>0.25\Longleftrightarrow
\beta>\frac{1}{4}.
\]
and all the $\beta>\frac{1}{4}$ can be achieved. In particular, when
$\beta=\frac{1}{3}$, where $\alpha^2=\frac{5}{18}=0.277777_{\dots}$,
one can find approximate value of $\tau\sim 6.73$ from numerical
result. In the picture, we have identified three special points:
$(1,0.5)$,$(3,1/3)$ and $(6.73, \frac{5}{18})$ which corresponds to
$\beta=1$, $\frac{1}{2}$ and $\frac{1}{3}$ respectively. The
corresponding graph of $f$ and $h=f_r$ for $\tau=6.73
(\beta=\frac{1}{3})$ is shown in figure \ref{P23D}.
%\begin{center}
%\includegraphics[height=4cm]{b3.eps}
%\includegraphics[height=4cm]{b3_2.eps}
%\end{center}
Finally, note that we are only interested when $\beta\le1$, or
equivalently when $\alpha^2\le 0.5$. However the picture suggests we
can even pass $\beta\le 1$ and solve for conic K\"{a}hler-Einstein
metric with cone angle $2\pi\beta>2\pi$ along the conic curve.

\section{Limit as $\beta$ goes to $1/4$}\label{limit}

\subsection{Metric Limit}
We know that $SU(2)$ acts on $\mathbb{P}^1$ naturally. As pointed
out in \cite{LiSu}, the following embedding is equivariant with
respect to the covering homomorphism $\phi: SU(2)\rightarrow
SO(3,\mathbb{R})$.
\begin{eqnarray*}
\Delta: \mathbb{P}^1&\longrightarrow &\mathbb{P}^2\\
{[U_0,U_1]}&\mapsto& [U_0^2+U_1^2,2iU_0U_1,i(U_0^2-U_1^2)].
\end{eqnarray*}
Here $SU(2)$ acts on $\mathbb{P}^2(1,1,4)$ by acting on the first
two variables:
\[
g\cdot [U_0,U_1,V]=[g\cdot (U_0,U_1), V].
\]
Note that
\[
\Delta(\mathbb{P}^1)=\{Z_1^2+Z_2^2+Z_3^2=0\}.
\]
Fix generators of $SU(2,\mathbb{C})$ to be standard Pauli matrices:
\[
Y_1=\left(\begin{array}{cc} 0&1\\-1&0\end{array}\right),\;
Y_2=\left(\begin{array}{cc} 0&i\\i&0\end{array}\right),\;
Y_3=\left(\begin{array}{cc} i&0\\0&-i\end{array}\right).
\]
Note that the commutator relation $[Y_1,Y_2]=2Y_3$ and cyclicly. So
by letting $\tilde{Y}_i=\frac{Y_i}{2}$, $\tilde{Y}_i$'s satisfy
$[\tilde{Y}_1,\tilde{Y}_2]=\tilde{Y}_3$ and cyclicly. For
simplicity, we will still use $\tilde{Y}_i$ to denote the vector
fields on $\mathbb{P}(1,1,4)$ corresponding to the infinitesimal
actions of $\tilde{Y}_i$. Then we have
\begin{lem}
When we restrict to $\mathbb{P}^1$, $\Delta_*\tilde{Y}_1=-T_u$,
$\Delta_*\tilde{Y}_2=T_v$, $\Delta_*\tilde{Y}_3=T_w$.
\end{lem}
\begin{proof}
$\Delta(1,0)=(1,0,i)=u+iv$ with $u=(1,0,0)$ and $v=(0,0,1)$. So
$w=u\times v=-(0,1,0)$.
\[
\Delta_*\tilde{Y}_i=(2(U_0\dot{U}_0+U_1\dot{U}_1),
2i(\dot{U}_0U_1+U_0\dot{U}_1), 2i(U_0\dot{U}_0-U_1\dot{U}_1))
\]
\begin{enumerate}
\item
$\tilde{Y}_1=\frac{1}{2}(U_1,-U_0)$, so
\[
\Delta_*\tilde{Y}_1=(0,-i(U_0^2-U_1^2),2iU_0U_1).
\]
In particular, $\tilde{Y}_1|_{(1,0)}=\frac{1}{2}(0,-1)$ and
$\Delta_*\tilde{Y}_1|_{(1,0)}=(0,-i,0)=ie_w$. So
$\Delta_*\tilde{Y}_1=-T_u$.
\item
$\tilde{Y}_2=\frac{1}{2}(iU_1,i U_0)$, so
\[
\Delta_*\tilde{Y}_2=(2iU_0U_1,-(U_0^2+U_1^2),0).
\]
In particular, $\tilde{Y}_2|_{(1,0)}=\frac{1}{2}(i,0)$,
$\Delta_*\tilde{Y}_2|_{(1,0,i)}=(0,-1,0)=e_w$. So
$\Delta_*\tilde{Y}_2=T_v$.
\item
$\tilde{Y}_3=\frac{1}{2}(i U_0,-iU_1)$, so
\[
\Delta_*\tilde{Y}_1=(i(U_0^2-U_1^2),0,-(U_0^2+U_1^2)).
\]
In particular, $\tilde{Y}_3|_{(1,0)}=\frac{i}{2}(1,0)$ and $
\Delta_*\tilde{Y}_3|_{(1,0,i)}=(i,0,-1)=-v+iu$. So
$\Delta_*\tilde{Y}_3=T_w$.\end{enumerate}
\end{proof}
We can define a function which classifies the $SU(2)$-orbits
\begin{eqnarray*}
\tilde{R}: \mathbb{P}(1,1,4)&\longrightarrow& [0,+\infty)\\
{[U_0,U_1,V]}&\mapsto &
\left(\frac{|U_0|^2+|U_1|^2}{|V|^{1/2}}\right)^{1/2}
\end{eqnarray*}
\begin{lem}
The generic orbit when $0<\tilde{R}<\infty$ is isomorphic to
$SU(2)/\mathbb{Z}_4\cong SO(3)/\mathbb{Z}_2$. The special orbit are
\[
{\rm Orb}_{\tilde{R}=0}={\rm Pt}=[0,0,1],\quad {\rm
Orb}_{\tilde{R}=\infty}=\mathbb{P}^1.
\]
\end{lem}
\begin{proof}
If $0<\tilde{R}<+\infty$, then $[U_0,U_1,V]$ is the same as
$[\sqrt{-1}^jU_0,\sqrt{-1}^jU_1,V]$, $j=1,2,3,4$. So the stabilizer
is isomorphic to $\mathbb{Z}_4$. The cases of special orbits are
clear.
\end{proof}
Now the $SU(2)$-invariant K\"{a}hler metric has the form
\[
g=dt^2+a^2 \sigma_1^2+b^2 \sigma_2^2+c^2 \sigma_3^2.
\]
Similar as the example \ref{P2FS} in Section \ref{orbitsP2}, we can
calculate the induced orbifold K\"{a}hler-Einstein metric by the
branch covering map:
\begin{eqnarray*}
\mathbb{P}(1,1,1)&\longrightarrow& \mathbb{P}(1,1,4)\\
{[Z_1,Z_2,Z_3]}&\mapsto& [Z_1,Z_2,Z_3^4].
\end{eqnarray*}
Because the metric is $SU(2)$ invariant, to write down the metric we
only need to calculate the length of the basic vector fields at the
the special point $(\tilde{R},0,1)$ in each $SU(2)$-orbit.
\begin{enumerate}
\item
$T_{\tilde{R}}|_{(\tilde{R},0,1)}=(1,0,0)$, $
|T_{\tilde{R}}|=\frac{1}{1+\tilde{R}^2}. $
\item
$ \tilde{Y}_1|_{(\tilde{R},0,1)}=\frac{1}{2}(0,-\tilde{R},0)$, $
a=|\tilde{Y}_1|_g=\frac{1}{2}\frac{\tilde{R}}{\sqrt{1+\tilde{R}^2}}$.
\item
$\tilde{Y}_2|_{(\tilde{R},0,1)}=\frac{1}{2}(0,i\tilde{R},0)$, $
b=|\tilde{Y}_2|_g=\frac{1}{2}\frac{\tilde{R}}{\sqrt{1+\tilde{R}^2}}$.
\item
$\tilde{Y}_3|_{(\tilde{R},0,1)}=\frac{1}{2}(i\tilde{R}, 0, 0)$, $
c=|\tilde{Y}_3|_g=\frac{1}{2}\frac{\tilde{R}}{1+\tilde{R}^2}$.
\end{enumerate}
Again, we can transform to the distance function:
\[
dt=-\frac{d\tilde{R}}{1+\tilde{R}^2}\; \&\; \tilde{R}(+\infty)=0
\Longrightarrow \tilde{R}=\tan (\pi/2-t), 0\le t\le \pi/2.
\]
By substituting $\tilde{R}$ into the expression of $a$, $b$ and $c$,
we get the data for $\mathbb{P}(1,1,4)$:
\[
a=b=\frac{1}{2}\sin\left(\frac{\pi}{2}-t\right)=\frac{1}{2}\cos(t).
\]
\[
c=\frac{1}{4}\sin (\pi-2t)=\frac{1}{4}\sin(2t).
\]
%\[
%\frac{dr}{dt}=\frac{1}{c}=\frac{4}{\sin (2t)}\Rightarrow
%e^{r/2}=\tan (t).
%\]
%\[
%f=ab=\frac{1}{4}\cos^2(t)=\frac{1}{4}\cdot \frac{1}{1+e^r}.
%\]
Note that in this case, $a/b\equiv 1$. This is very different from
the case where $\beta>1/4$. For the latter, $a< b$ except on the
special fibre ${\rm Orb}_{R=1}\cong\mathbb{P}^1$ where $a=b$.
Moreover, the boundary condition now becomes
\begin{eqnarray*}
a(t)=b(t)&=& 1/2+O(t^2)\\
 c(t)&=& \frac{1}{2}t+O(t^3)
\end{eqnarray*}
On the other end where $t_*=\pi/2$, $a(\pi/2)=b(\pi/2)=c(\pi/2)=0$.
Geometrically, the special fibre ${\rm Orb}_{R=0}\cong\mathbb{RP}^2$
shrinks to a point as $\beta\rightarrow 1/4$. If we do the same
transformation that $dr/dt=1/c$, the range of $r$ will becomes
$(-\infty, +\infty)$ instead of $(-\infty,0)$ because $c(t_*)=0$.

Next we give the numerical results which show that the metric
$\omega_\beta$ converges to the orbifold K\"{a}hler-Einstein metric
on $\mathbb{P}(1,1,4)$.

First we integrate the identity $dr/dt=1/c$ numerically and plot the
relation between the boundary value $b(t_*)^2=1/\tau$ and
$t_{\max}=t_*$. We see that the maximal value for $t$ is an
increasing function of $\tau$. As $\tau\rightarrow +\infty$, or
equivalently as $\beta\rightarrow 1/4$, $t_{max}=t_*$ converges to
$\pi/2$.
\begin{figure}[h]
  \begin{center}
    \subfigure[short range]{\label{tautmax-a}\includegraphics[height=4cm]{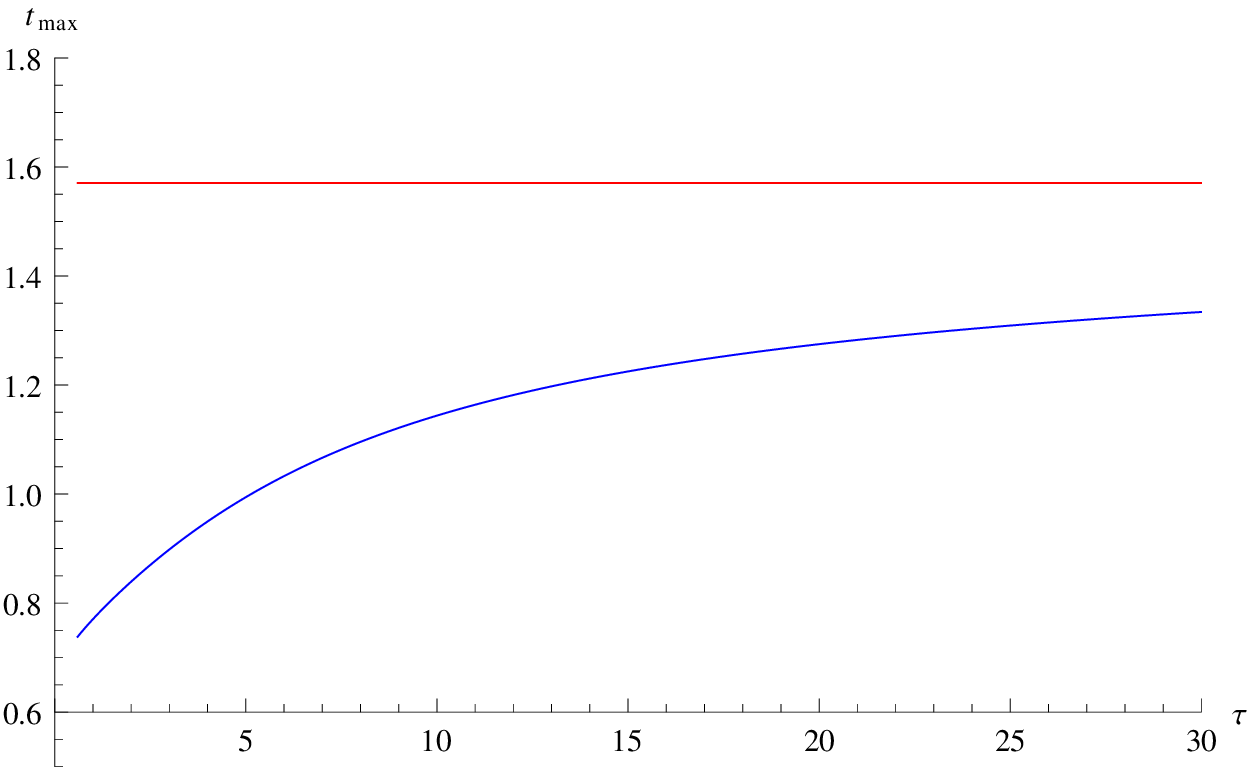}}
    \subfigure[long range]{\label{tautmax-b}\includegraphics[height=4cm]{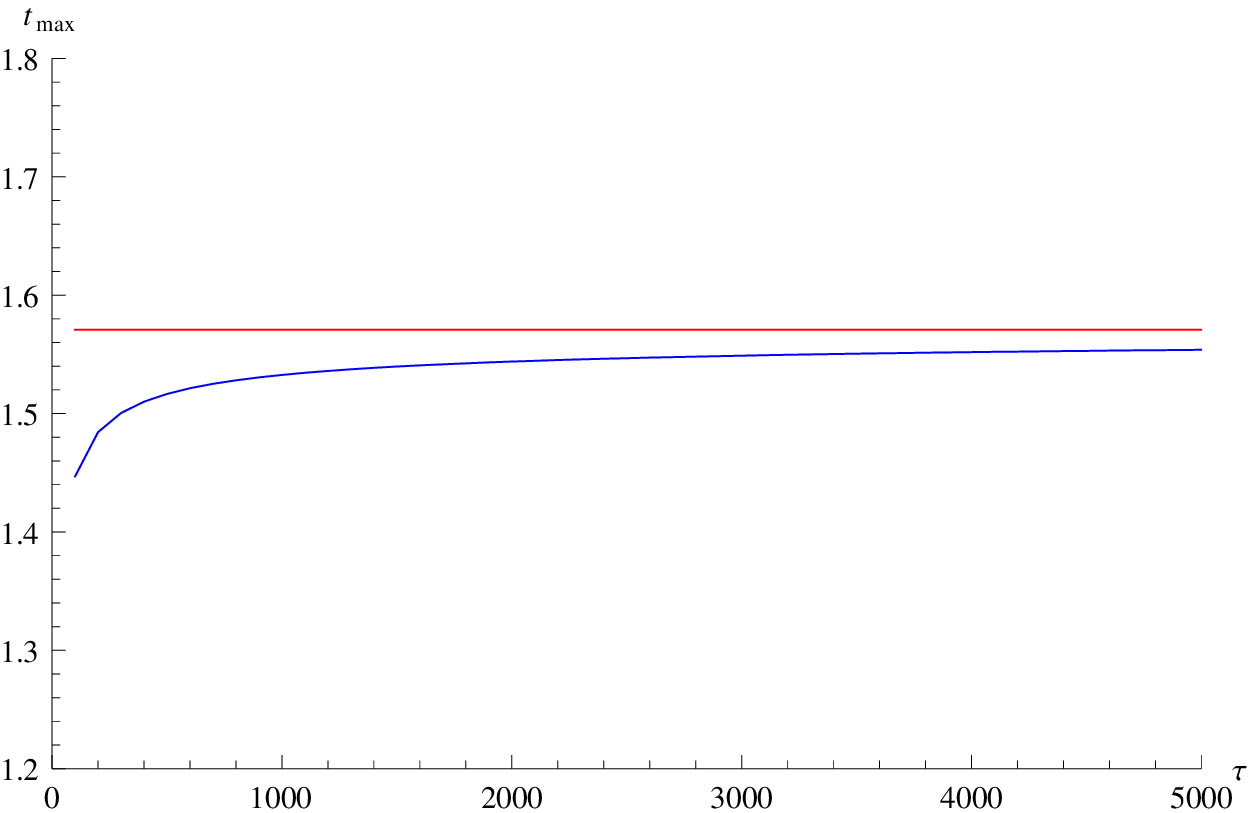}}
  \end{center}
  \caption{$(\tau, t_{max}=t_*)$}
  \label{tautmax}
\end{figure}
%\begin{center}
%\includegraphics[height=4cm]{short.eps}
%\includegraphics[height=4cm]{long.eps}
%\end{center}
Note that the coordinate $t$ is the distance function from the
special orbit $\mathbb{P}^1$. So $t$ is a geometrically meaningful
coordinate in contrast with $r$ which is only an auxiliary
coordinate. So we can get a good convergence when we look the
data as functions $t$.

\begin{figure}[h]
  \begin{center}
    \subfigure[$(t,f=ab)$]{\label{convdata-a}\includegraphics[height=4cm]{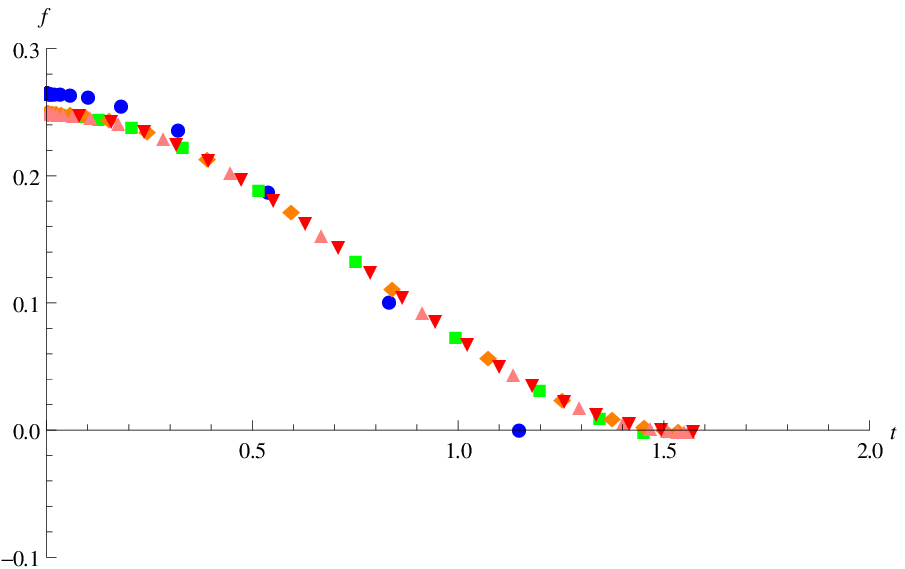}}
    \subfigure[$(t,a/b)$]{\label{convdata-b}\includegraphics[height=4cm]{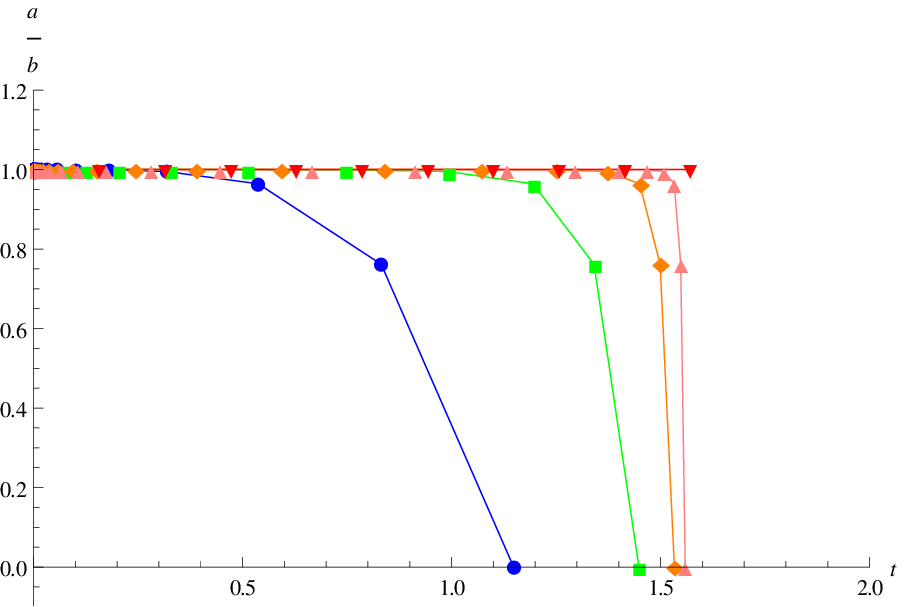}}  \\
    \subfigure[$(t,c^2)$]{\label{convdata-c}\includegraphics[height=4cm]{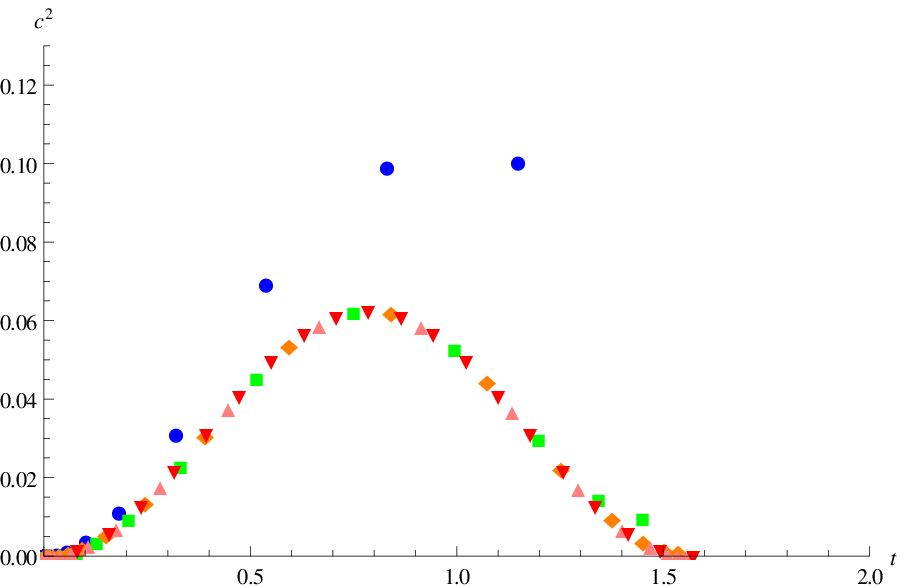}}
  \end{center}
  \caption{Convergence of data}
  \label{convdata}
\end{figure}

Now we can plot the graph of the data set $(f=ab,R=a/b,-f_r=c^2)$ as
the function of $t$ instead of $r$. (See \eqref{2sets}). 
Figure \ref{convdata} shows the data for four $\tau$'s: $\tau=10^i$ for
$i=1,2,3,4$. The corresponding colors and markers are ``Blue Round",
``Green Square", ``Orange Diamond", ``Pink Triangle" for $i=1,2,3,4$
respectively. The ``Red Upside-down Triangle" represent the data for
$\mathbb{P}(1,1,4)$ where
\[
f(t)=a(t)b(t)=\frac{1}{4}\cos^2(t),\quad R=\frac{a}{b}\equiv 1,
\quad c^2(t)=\frac{1}{4}\sin^2(2t).
\]
%\begin{center}
%\includegraphics[height=4cm]{limitf.eps}
%\includegraphics[height=4cm]{limitR.eps}
%\includegraphics[height=4cm]{limitc.eps}
%\end{center}
One can see that the data for $\tau$ large fits with the data for
$\mathbb{P}(1,1,4)$ very well. Again, we know that $\tau$ going to
$+\infty$ is equivalent to $\beta$ going to $1/4$. So the numerical
result implies the expected result: as $\beta\rightarrow 1/4$, the
metric $\omega_\beta$ converges to the orbifold K\"{a}hler-Einstein
metric $\hat{\omega}_{KE}$ on $P(1,1,4)$.

\subsection{$\mathbb{Z}_2$-quotient of Eguchi-Hanson as the Bubble}

\begin{figure}[h]
  \begin{center}
    \subfigure[$(r,f)$]{\label{bubble-a}\includegraphics[height=4cm]{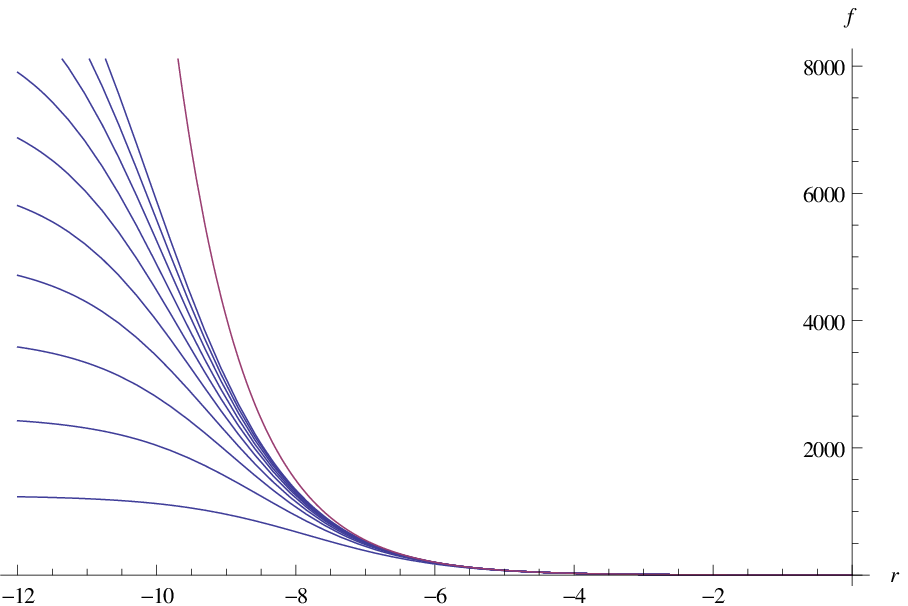}}
    \subfigure[$(r,f_r)$]{\label{bubble-b}\includegraphics[height=4cm]{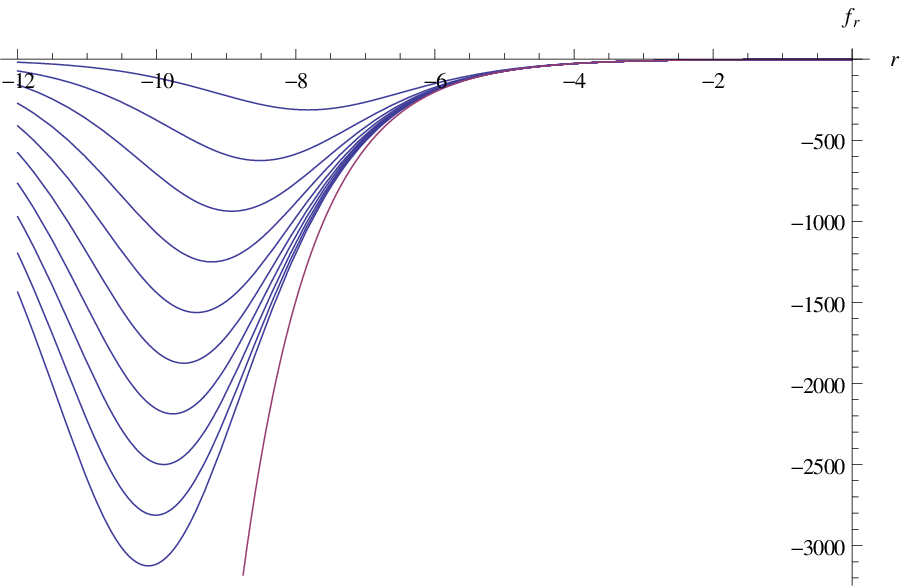}} 
  \end{center}
  \caption{Bubbling}
  \label{bubble}
\end{figure}

As pointed out by Dr. H-J. Hein and Professor Lebrun, if we rescale the metric near the orbit ${\rm Orb}_{R=0}=\mathbb{RP}^2$ appropriately, then the rescaled metrics should converge to another
well known metric which is the $\mathbb{Z}_2$ quotient of the Eguchi-Hanson metric. This kind of metrics was studied in 
much generality by Stenzel \cite{Sten}. It's easy to see this convergence from the discussion in Section \ref{secKEeq} and the following numerical results. For this we use the 
explicit description of this metric in \cite[Section 7]{Sten}, which says that, away from the $\mathbb{RP}^2$ the $\mathbb{Z}_2$-quotient of Eguchi-Hanson metric can be pulled back to an $SO(3)$ 
invariant metric on $(0,\infty)\times SO(3)$ with the following expression:
\begin{equation}\label{EHmetric}
g=\cosh s (ds)^2+\sinh s \tanh s (X_1^*)^2+\cosh s ((X_2^*)^2+(X_3^*)^2).
\end{equation}
As before, we can let $a^*(s)=\sqrt{\sinh s \tanh s}$, $b^*(s)=c^*(s)=\sqrt{\cosh s}$. Let $t^*$ be the distance function to the orbit $\mathbb{RP}^2$. Then from
\eqref{EHmetric}, we see the following relation:
\[
\frac{ds}{dt^*}=\frac{1}{\sqrt{\cosh s}}=\frac{1}{c^*(s)},\quad \frac{a^*}{b^*}(s)=\tanh s.
\]
If we compare these identities with \eqref{coordtr} and \eqref{aoverb}, we see that the coordinate $r$ is preserved under this convergence. 
In other words, $r=-s$ and $\frac{a^*}{b^*}=\frac{a}{b}=-\tanh r$. To prove the convergence,
we only need to prove the convergence of rescaled data as functions of $r$. Note that, since the length scale of ${\rm Orb}_{R=0}=\mathbb{RP}^2$ is $1/\sqrt{\tau}$ as 
$\tau\rightarrow+\infty$ (equivalently as $\beta\rightarrow 1/4$), we need to use the scale factor $\tau$ to rescale the metric back. So we need to show the following convergence.
\[
\lim_{\tau\rightarrow +\infty}f\tau=\lim_{\tau\rightarrow +\infty} a(r,\tau)b(r,\tau)\cdot \tau=a^*b^*=-\sinh r,
\]
\[
\lim_{\tau\rightarrow +\infty} f_r \tau=\lim_{\tau\rightarrow+\infty}-c(r,\tau)^2\tau=-(c^*)^2=-\cosh r.
\]
Figure \ref{bubble} shows the convergence of numerical data for $\tau=5000*i$, $i=1,\cdots, 10$.

\section{Data of $\bP^1\times\bP^1$ and associated Sasaki-Einstein
metric}\label{p1p1}

We have the following Segre embedding of $\bP^1\times\bP^1$ into
$\bP^3$ by the complete linear system $|H_1+H_2|$ where $H_1$ and
$H_2$ are the hyperplane divisors of the two factors of $\bP^1$
respectively.
\begin{eqnarray}\label{segre}
\phi: \mathbb{P}^1\times\mathbb{P}^1&\longrightarrow& \bP^3\\
([U_0,U_1],[V_0,V_1])&\mapsto& [U_0V_0+U_1V_1, \cI(U_0V_0-U_1V_1),
U_0V_1+U_1V_0,U_0V_1-U_1V_0].\nonumber
\end{eqnarray}
Note that
\[
\phi(\bP^1\times\bP^1)=\{[Z_1,Z_2,Z_3,Z_4]\in\bP^3;
Z_1^2+Z_2^2+Z_3^2=Z_4^2\}.
\]
\begin{lem}
Let $p_i: \bP^1\times\bP^1\rightarrow \bP^1$ be the projection to the $i-th$ $\bP^1$-factor and $\omega_{\bP^N}$ denote the standard Fubini-Study metric on
$\bP^N$ in the cohomology class $2\pi c_1(\mathcal{O}_{\bP^N}(1))$,
then the Segre embedding $\phi$ satisfies
\[
\phi^*\omega_{\bP^2}=p_1^*\omega_{\bP^1}+p_2^*\omega_{\mathbb{P}^1}=:\tilde{\omega}.
\]
\end{lem}
\begin{proof}
 This follows from the following formula:
 \begin{eqnarray*}
 &&p_1^*\omega_{\bP^1}+p_2^*\omega_{\bP^1}=\sddbar\log((|U_0|^2+|U_1|^2)(|V_0|^2+|V_1|^2))\\
 &=&\sddbar\log\left(|U_0V_0+U_1V_1|^2+|\cI(U_0V_0-U_1V_1)|^2+|U_0V_1+U_1V_0|^2+|U_0V_1-U_1V_0|^2\right).
 \end{eqnarray*}
\end{proof}
Now $SO(3)$ acts on $\mathbb{C}^4$ by
\[
g\cdot (Z_1,Z_2,Z_3,Z_4)=(g\cdot (Z_1,Z_2,Z_3), Z_4).
\]
This induces an action of $SO(3)$ on $\phi(\bP^1\times\bP^1)$.

We will calculate the data associated with the product metric
$\tilde{\omega}:=p_1^*\omega_{\bP^1}+p_2^*\omega_{\bP^1}$. Use the
similar method as in Section \ref{orbitsP2} we use the following
notation:
\[
(Z_1,Z_2,Z_3,Z_4)\sim e^{-i{\rm
Arg}(Z_1^2+Z_2^2+Z_3^2)/2}(Z_1,Z_2,Z_3,Z_4)=:(u+iv, z_4).
\]
Here $u,v\in \mathbb{R}^3$, $z_4\in \mathbb{C}$. In this notation,
we have
\[
\phi(\bP^1\times\bP^1)=\{(u+iv,z_4); |u|^2-|v|^2=z_4^2, 0\neq
(u+iv,z_4)\in \mathbb{C}^3\times\mathbb{R}\}/\mathbb{R}^{\times}.
\]
We can calculate the infinitesimal vector field of basis of $so(3)$,
at point $(u+iv, \sqrt{|u|^2-|v|^2})$:
\[
T_u=(-\sqrt{-1}|v|e_w,0),\quad T_v=(|u|e_w,0),\quad T_w=(-|u|e_v+\sqrt{-1}
|v|e_u,0).
\]
As in Section \ref{P2FS}, we define $R=\frac{|v|}{|u|}$ and
calculate the radial vector field as
\[
\quad T_R=\left(\cI |u|e_v,-\frac{|u|R}{\sqrt{1-R^2}}\right). \] 
Here for clarify, we will use $\overline{T}_u$, $\overline{T}_v$, $\overline{T}_w$ and
$\overline{T}_R$ to denote the tangent vector in
$T_{[u+iv]}\mathbb{P}^2$ determined by $T_u, T_v, T_w, T_R$
respectively.
The lengths of these tangent vectors in
$T_{[u+iv,iz_3]}\phi(\bP^1\times\bP^1)$ can be calculated as in
Example \ref{P2FS}:
\[
|\overline{T}_u|_{\tilde{\omega}}^2=\frac{R^2}{2},\quad
|\overline{T}_v|_{\tilde{\omega}}^2=\frac{1}{2}, \quad
|\overline{T}_w|_{\tilde{\omega}}^2=\frac{1-R^2}{2},\quad
|\overline{T}_R|_{\tilde{\omega}}^2=\frac{1}{2(1-R^2)}.
\]
By transforming the variable $R$ into the distance variable
$\tilde{t}$ under the metric $\tilde{\omega}$, we get:
\begin{equation}\label{tfuncR}
\frac{d\tilde{t}}{dR}=\frac{1}{2(1-R^2)}\;\&\;
\tilde{t}(1)=0\Longrightarrow R(\tilde{t})=\cos(\sqrt{2}\tilde{t}),
0\le \tilde{t}\le \frac{\pi}{2\sqrt{2}}.
\end{equation}
\[
|\overline{T}_u|_{\tilde{\omega}}=\frac{1}{\sqrt{2}}\cos(\sqrt{2}\tilde{t}),\quad
|\overline{T}_v|_{\tilde{\omega}}=\frac{1}{\sqrt{2}},\quad
|\overline{T}_w|_{\tilde{\omega}}=\frac{1}{\sqrt{2}}\sin(\sqrt{2}\tilde{t}).
\]
Note that $\tilde{\omega}=p_1^*\omega_{\bP^1}+p_2^*\omega_{\bP^1}$ has Ricci
curvature equal to 4. To normalize Ricci curvature to be 6, we just
need to rescale the metric. So by letting
$\omega=\frac{2}{3}\omega_{\bP^1\times\bP^1}$ and redefining
$t=\sqrt{2}\tilde{t}/\sqrt{3}$ we get the following result, which are
the same data as in Example \ref{P1P1FS}
\begin{equation}
a=\frac{1}{\sqrt{3}}\cos(\sqrt{3}t),\quad b=\frac{1}{\sqrt{3}},\quad
c=\frac{1}{\sqrt{3}}\sin(\sqrt{3}t);\quad 0\le t\le
\frac{\pi}{2\sqrt{3}}.
\end{equation}
Let ${\rm Aff}(\bP^1\times\bP^1)$ be the affine cone over
$\phi(\bP^1\times\bP^1)\subset\bP^3$:
\[
{\rm Aff}(\bP^1\times\bP^1)=\{(Z_1,Z_2,Z_3,Z_4)\in \mathbb{C}^4;
Z_1^2+Z_2^2+Z_3^2=Z_4^2\}.
\]
In the following, we use $L$ to denote the total space of the line
bundle $p_1^*\mathcal{O}(-H_1)+p_2^*\mathcal{O}(-H_2)$. Then
$L=Bl_{0}{\rm Aff}(\bP^1\times\bP^1)$. In other words, the zero
section $S_0$ of $L$ can be blow-down to get a singular variety
$L/S_0$ which is isomorphic to the affine cone ${\rm
Aff}(\bP^1\times\bP^1)$. Moreover, line bundle $L$ has a Hermitian
metric $h:=h_{\bP^1\times\bP^1}$ whose curvature is
$-\tilde{\omega}=-(p_1^*\omega_{\bP^1}+p_2^*\omega_{\bP^1})$, i.e.
we have the identity:
\begin{equation}
-\sddbar\log
h=-\tilde{\omega}=-(p_1^*\omega_{\bP^1}+p_2^*\omega_{\bP^1}).
\end{equation}
Now $h: L\ni s\rightarrow |s|_{h}^2$ is a smooth function on $L$
which induces a smooth function $h$ on $L/S_0\cong {\rm
Aff}(\bP^1\times\bP^1)$. Up to a scaling factor, we see that
\begin{eqnarray*}
h: {\rm Aff}(\bP^1\times\bP^1)&\rightarrow&\mathbb{R}^{\ge 0}\\
(Z_1,Z_2,Z_3,Z_4)&\mapsto & |Z|^2=|Z_1|^2+|Z_2|^2+|Z_3|^2+|Z_4|^2.
\end{eqnarray*}
Define $M^5\subset L$ to be the unit circle bundle, i.e. $M^5=\{s\in
L; |s|_h^2=1\}$. Then
\[
M^5\cong \{(Z_1,Z_2,Z_3,Z_4); Z_1^2+Z_2^2+Z_3^2=Z_4^2,|Z|^2=1\}={\rm
Aff}(\bP^1\times\bP^1)\cap S^7 .
\]
We know that there exists a Sasaki-Einstein metric on $M^5$. Now we will calculate this Sasaki-Einstein metric on $M^5$ by calculating
the data in the sense of \cite{Conti}. To do this we will first
calculate the metric on $M^5$ induced by the standard Euclidean
metric on $\mathbb{C}^4$. Then we modify the metric appropriately
(rescale it in different directions) to get the desired
Sasaki-Einstein metric.
\begin{lem}
On $M^5=S^7\cap {\rm Aff}(\bP^1\times\bP^1)$, we have
\begin{equation}\label{lengthuv}
|u|=\frac{1}{\sqrt{2}},\quad
|v|=\frac{\cos(\sqrt{2}\tilde{t})}{\sqrt{2}}.
\end{equation}
\end{lem}
\begin{proof}
On $M^5$, we have the identities $|u|^2-|v|^2=z_4^2$ and
$|u|^2+|v|^2+|z_4|^2=1$. So we get $|u|=1/\sqrt{2}$. The second
identity follows from \eqref{tfuncR} and $|v|=R |u|$.
\end{proof}
Now $G=SO(3)\times U(1)$ acts on $\mathbb{C}^4$ by
\[
(g,e^{i\theta})\cdot (Z_1,Z_2,Z_3,Z_4)=\left(e^{i\theta} g(Z_1,Z_2,Z_3), e^{i\theta}Z_4\right).
\]
The generic orbit is of codimension 1. ${\rm Aff}(\bP^1\times\bP^1)$ is $G$-invariant under this action. 
Fix the standard basis of $so(3)\oplus u(1)$ by adjoining the
generator $X_4$ of $u(1)$ to the standard basis of $so(3)$ used
above. We will denote the infinitesimal vector fields by the same
notation. So we have
\[
X_1=T_u,X_2=T_v, X_3=T_w, X_4=(-v+iu,iz_4).
\]
\begin{prop}\label{unscaled}
Considering $M^5$ as a submanifold in $(\mathbb{C}^4, g_{{\rm flat}})$,
$T_{(u+iv,z_4)}\mathbb{C}^4=\mathbb{R}^8$ has an orthonormal basis
given by
 \begin{eqnarray*}
 \partial_{\tilde{r}}&=&(u+iv,z_4)\\
\tmfe_0&=&\partial_\theta=X_4=(-v+iu,i z_4),\\
\tilde{\mathfrak{e}}_1&=&\partial_{\tilde{t}}=-\sqrt{2}\sin(\sqrt{2}\tilde{t})T_R=(-\cI\sin(\sqrt{2}\tilde{t})e_v,\cos(\sqrt{2}\tilde{t})),\\
\tilde{\mathfrak{e}}_2&=&\frac{\sqrt{2}}{\sin(\sqrt{2}\tilde{t})}(-X_3+\cos(\sqrt{2}\tilde{t})X_4)=(\sin(\sqrt{2}\tilde{t})
e_v,
\cI\cos(\sqrt{2}\tilde{t})),\\
\tilde{\mathfrak{e}}_3&=&\frac{\sqrt{2}}{\cos(\sqrt{2}\tilde{t})}X_1=(-\cI e_w,0),\\
\tilde{\mathfrak{e}}_4&=&\sqrt{2}X_2=(e_w, 0).
\end{eqnarray*}
We have the relation
\begin{equation*}
J\partial_{\tilde{r}}=i\partial_{\tilde{r}}=\partial_\theta, \quad J
\tmfe_1=i\tmfe_1=\tmfe_2,\quad J\tmfe_3=i\tmfe_3=\tmfe_4.
\end{equation*}
Under the induced metric on $M^5$ by the standard Euclidean metric
on $\mathbb{C}^4$, $T_{(u+iv,z_4)}M^5$ has an orthonormal basis
$\{\partial_\theta, \tmfe_1,\tmfe_2,\tmfe_3, \tmfe_4\}$. Moreover,
let $S^1\rightarrow M^5\rightarrow \bP^1\times\bP^1$ be the
fibration structure. Then the vertical unit vector field is
generated by $\partial_\theta$, and the space of horizontal vector
fields in the tangent space has an orthonormal basis consisting of
$\{\tmfe_1,\tmfe_2,\tmfe_3,\tmfe_4\}$.
\end{prop}
\begin{proof}
First it's easy to see that $J\partial_{\tilde{r}}=\partial_\theta$
and $\partial_{\tilde{r}}\perp
Span(\{\partial_\theta,X_i,i=1,2,3,4\})$. We can also verify that
$Span\{X_1,X_2\}\perp Span(\{\partial_\theta, X_3,X_4\})$ and
$T_R\perp Span(\{\partial_\theta, X_i,i=1,2,3,4 \})$. The Lemma
follows by orthonormalization.
\end{proof}

\begin{lem}
Considering $h$ as a smooth function on $L$ as above, Sasaki-Einstein metric on $M^5$ is given by
\[
g_{SE}=\frac{1}{2}\left.(\sddbar h^{2/3})(\cdot, J\cdot)\right|_{M^5}.
\]
\end{lem}
\begin{proof}
If $M^5$ is a Sasaki-Einstein metric, then the metric cone $C(M^5)$
is a Ricci-flat K\"{a}hler metric. In our case, $C(M^5)\cong L/S_0$ as
the affine variety with an isolated singular point. So we only need to construct the rotationally
symmetric Ricci-flat K\"{a}hler metric on $C(M^5)\cong L/S_0$ and restrict to
$M^5\cong \{h=1\}\cap C(M^5)$ to get the Sasaki-Einstein metric on
$M^5$.

In general, assume $L\rightarrow D_0$ be a line bundle with a
Hermitian metric $h$ such that $\sddbar\log h=\tilde{\omega}$ is a
K\"{a}hler-Einstein metric, satisfying
$Ric(\tilde{\omega})=\tau\tilde{\omega}$. Then we can define the
rotationally symmetric K\"{a}hler metric on the total space on $L/S_0$
using the potential $h^{\delta}$, i.e. we define
\begin{equation}\label{expdelta}
\Omega_\delta=\sddbar h^{\delta}=\delta
h^{\delta}\tilde{\omega}+\delta^2 h^{\delta}
\frac{\nabla\xi\wedge\overline{\nabla\xi}}{|\xi|^2}.
\end{equation}
The Ricci curvature of $\Omega_\delta$ on $L\backslash D_0$ is equal
to
\begin{eqnarray*}
Ric(\Omega_\delta)&=&-\sddbar\log\Omega_\delta^n=-(d+1)\sddbar \log
h^\delta+Ric(\tilde{\omega})\\
&=&\pi^*(-(d+1)\delta\tilde{\omega}+\tau\tilde{\omega}).
\end{eqnarray*}
This is zero if and only if $\delta=\tau/(d+1)$. In our case,
$\tau=2$, $d=2$. So $\delta=2/3$.
\end{proof}
\begin{thm}
The Sasaki-Einstein metric on $M^5$ has an orthonormal basis given
by
 \begin{eqnarray*}
\mathfrak{e}_0&=&\frac{3}{2}X_4,\\
\mathfrak{e}_1&=&\partial_{t}=-\sqrt{3}\sin(\sqrt{3}t)T_R,\\
\mathfrak{e}_2&=&\frac{\sqrt{3}}{\sin(\sqrt{3}t)}(-X_3+\cos(\sqrt{3}t)X_4),\\
\mathfrak{e}_3&=&\frac{\sqrt{3}}{\cos(\sqrt{3}t)}X_1,\\
\mathfrak{e}_4&=&\sqrt{3}X_2.
\end{eqnarray*}
\end{thm}
\begin{proof}
First note that, the induced metric on $M^5$ by flat metric is given
by $\frac{1}{2}\sddbar h$. By the formula \eqref{expdelta}, we see that if we
change the potential from $h$ to $h^{\delta}$, then the vertical
metric scales by $\delta^2$, and the horizontal part of the metric
scales by $\delta$. Since $\delta=2/3$ now, the Theorem follows from
Proposition \ref{unscaled}.
\end{proof}
\begin{cor}\label{P1P1SU}
\begin{enumerate}
\item
Under the Sasaki-Einstein metric on $M^5$, there is an orthonormal basis of $T^*M^5$ given by
\begin{eqnarray*}
&&\alpha:=\frac{2}{3}(X_4^*+\cos(\sqrt{3}t)X_3^*), \\
&&\mfe^1=dt,\quad \mfe^2=-\frac{\sin(\sqrt{3}t)}{\sqrt{3}}X_3^*,\\
&&\mfe^3=\frac{\cos(\sqrt{3}t)}{\sqrt{3}}X_1^*,\quad
\mfe^4=\frac{1}{\sqrt{3}}X_2^*.
\end{eqnarray*}
\item
If we define $\omega_1=e^1\wedge e^2+e^3\wedge e^4$,
$\omega_2=e^1\wedge e^3+e^4\wedge e^2$ and $\omega_3=e^1\wedge
e^4+e^2\wedge e^3$, the following identities hold:
\[
d\alpha=2\omega_1,\quad d\omega_2=-3\alpha\wedge
\omega_3+2X_4^*\wedge \omega_3,\quad
d\omega_3=3\alpha\wedge\omega_2-2X_4^*\wedge\omega_2.
\]
This gives the $SU(2)$ structure in the sense of \cite{Conti}.
\end{enumerate}
\end{cor}
\begin{rem}
 The item 2 in Corollary \ref{P1P1SU} follows from Item 1 and the fomula $dX_i^*=-\epsilon_{ijk}X_j^*\wedge X_k^*$. As explained in \cite{Conti}, because we are using the $G$-invariant
 forms on $G\times (t_{-},t_{+})$ to represent the data, there is an extra term $2 X_4^*\wedge \omega_3$. The coefficient $2$ comes from 
 the fact that $(e^{i\theta})^*\mathcal{S}=e^{2i\theta}\cdot\mathcal{S}$ where we use $\mathcal{S}$ to denote the nonwhere vanishing holomorphic volume form on $\mathcal{M}={\rm Aff}(\bP^1\times\bP^1)$
 which can be given by the Poincar\'{e} residue formula:
 \[
 \mathcal{S}=Res_{\mathcal{M}}(dZ_1\wedge dZ_2\wedge dZ_3\wedge dZ_4)=-\frac{dZ_1\wedge dZ_2\wedge dZ_3}{Z_4}.
 \]
\end{rem}
\begin{rem}
By the similar calculation, we can calculate the data associated
on the standard round $S^5$ under the $SO(3)$ action:
\[
g\cdot (Z_1,Z_2,Z_3)=\left(g (Z_1,Z_2,Z_3)\right).
\]
The result is as follows. For the orthonormal basis of $TS^5$, we
have
\begin{eqnarray*}
&&\mfe_0=\partial_\theta=X_4,\\
&&\mfe_1=\partial_t,\quad \mfe_2=\frac{1}{\sin(2t)}(-X_3+\cos(2t)X_4),\\
&&\mfe_3=\frac{X_1}{\sin\left(\frac{\pi}{4}-t\right)},\quad
\mfe_4=\frac{X_2}{\cos\left(\frac{\pi}{4}-t\right)}.
\end{eqnarray*}
So the corresponding orthonormal basis of $T^*S^5$ is
\begin{eqnarray*}
&&\alpha:=\mfe^0=X_4^*+\cos(2t)X_3^*, \\
&&\mfe^1=dt,\quad \mfe^2=-\sin(2t)X_3^*,\\
&&\mfe^3=\sin\left(\frac{\pi}{4}-t\right)X_1^*,\quad
\mfe^4=\cos\left(\frac{\pi}{4}-t\right)X_2^*.
\end{eqnarray*}
The corresponding $SU(2)$-structural equations are:
\[
d\alpha=2\omega_1,\quad
d\omega_2=-3\alpha\wedge\omega_3+3X_4^*\wedge\omega_3,\quad
d\omega_3=3\alpha\wedge\omega_2-3X_4^*\wedge\omega_2.
\]
\end{rem}
\begin{rem}
There is a statement in Theorem 1 in \cite{Conti}: ``There is no
solution of (23) that defines an Einstein-Sasaki metric on a compact
manifold". The above two special examples show that this statement
is wrong. By going through the proof, we find that the error happens
in Lemma 4, where, in the second case, the assumption $q\neq 0$ is
made. In our notation, this implies the isotopy group of special
orbit has a generator whose $X_4$-component is nonzero. But this is
not true in the above examples. Actually, it's easy to verify that
\begin{enumerate}
 \item For $t=0$, $H_{-}\cong U(1)$ with Lie algebra $\mathfrak{h}=\langle -X_3+X_4 \rangle$.

 \item For $t=\frac{\pi}{2\sqrt{3}}$, $H_{+}\cong U(1)=U(1)_1$ with Lie algbra $\mathfrak{h}=\langle X_1 \rangle$. 
\end{enumerate}
Because the action $U(1)_1$ has generator $X_1$ which has no contribution from
$X_4$, so $q=0$ for $H_{+}$. It would be interesting to classify the missing cohomogeneity
one Sasaki-Einstein 5-manifolds for which $q=0$.
\end{rem}

\section{Appendix}

The following are the codes of Mathematica generating the figures appeared above.

\begin{enumerate}
 \item Figure \ref{P2}

  s = NDSolve[ \{f'[t] == h[t], h'[t] == 12 f[t]*h[t] + 2 Coth[2 t]*h[t] - h[t]\^{}2/f[t], 
     f[-10\^{}(-5)] == 10\^{}(-5), h[-10\^{}(-5)] == -1\}, \{f, h\}, \{t, -100, -10\^{}(-5)\}]; 
   
Plot[Evaluate[h[t] /. s], \{t, -3, -10\^{}(-5)\}, AxesLabel -$>$ \{r, Subscript[f, r]\}, 
 PlotRange -$>$ \{\{-3, 0\}, \{-1.2, 0.2\}\}, PlotLabel -$>$ P\^{}2] 
 
 Plot[Evaluate[h[t] /. s], \{t, -3, -10\^{}(-5)\}, AxesLabel -$>$ \{r, Subscript[f, r]\}, PlotRange -$>$ \{\{-3, 0\}, \{-1.2, 0.2\}\}, PlotLabel -$>$ P\^{}2]
 
 \item Figure \ref{P1P1}
 
 s1 = NDSolve[\{f'[t] == h[t], h'[t] == 12 f[t]*h[t] + 2 Coth[2 t]*h[t] - h[t]\^{}2/f[t], 
    f[-10\^{}(-5)] == 10\^{}(-5)/3, h[-10\^{}(-5)] == -1/3\}, \{f, 
    h\}, \{t, -100, -10\^{}(-5)\}];
    
    Plot[Evaluate[f[t] /. s1], \{t, -10, -10\^{}(-5)\}, 
 PlotRange -$>$ \{\{-10, -10\^{}(-5)\}, \{-0.1, 0.4\}\}, 
 PlotLabel -$>$ Superscript[P, 1]$\times$ Superscript[P, 1], 
 AxesLabel -$>$ \{r, f\}]
 
 Plot[Evaluate[h[t] /. s1], \{t, -5, -10\^{}(-5)\}, 
 AxesLabel -$>$ \{r, Subscript[f, r]\}, 
 PlotRange -$>$ \{\{-5, 0\}, \{-0.4, 0.1\}\}, 
 PlotLabel -$>$ Superscript[P, 1]$\times$Superscript[P, 1]]
 
 \item Figure \ref{taualpha}
 \begin{enumerate}
 \item Figure \ref{taualpha-a}
 
 Array[p,300];
 For[i = 0, i $<$ 300, 
  i++, \{v = NDSolve[\{f'[t] == h[t], 
      h'[t] == 12 f[t]*h[t] + 2 Coth[2 t]*h[t] - h[t]\^{}2/f[t], 
      f[-10\^{}(-5)] == 10\^{}(-5)/(0.5 + 0.1*i), 
      h[-10\^{}(-5)] == -1/(0.5 + 0.1*i)\}, \{f, h\}, \{t, -500, -10\^{}(-5)\}], 
   p[i + 1] = Evaluate[f[-500] /. v]\}];
 
 ListLinePlot[\{Table[\{0.5 + 0.1*i, Extract[p[i + 1], 1]\}, \{i, 0, 300\}],
   Table[\{0 + 0.1*i, 0.25\}, \{i, 0, 300\}], 
  Table[\{1, 0.05*i\}, \{i, 0, 10\}], Table[\{3, 1/3*0.1*i\}, \{i, 0, 10\}], 
  Table[\{0.1*i, 0.5\}, \{i, 0, 10\}], Table[\{0.3*i, 1/3\}, \{i, 0, 10\}], 
  Table[\{0.673*i, 5/18\}, \{i, 0, 10\}], 
  Table[\{6.73, 5/180*i\}, \{i, 0, 10\}]\}, 
 PlotRange -$>$ \{\{0, 22\}, \{0.2, 0.8\}\}, 
 AxesLabel -$>$ \{$\backslash$[Tau], $\backslash$[Alpha]\^{}2\}]
 
 \item Figure \ref{taualpha-b}
 
 Array[p, 50];
For[i = 1, i $<$ 51, 
  i++, \{v = 
    NDSolve[\{f'[t] == h[t], 
      h'[t] == 12 f[t]*h[t] + 2 Coth[2 t]*h[t] - h[t]\^{}2/f[t], 
      f[-10\^{}(-5)] == 10\^{}(-5)/(100*i), h[-10\^{}(-5)] == -1/(100*i)\}, \{f, 
      h\}, \{t, -500, 10\^{}(-4)\}], p[i] = Evaluate[f[-500] /. v]\}];
 
 ListLinePlot[\{Table[\{100*i, Extract[p[i], 1]\}, \{i, 1, 50\}], 
  Table[\{100*i, 1/4\}, \{i, 1, 50\}]\}, 
 PlotRange -$>$ \{\{0, 5000\}, \{0.249998, 0.250008\}\}, 
 AxesLabel -$>$ \{$\backslash$[Tau], $\backslash$[Alpha]\^{}2\}] 
 
 \end{enumerate}
 
 \item Figure \ref{P23D}
 
 v = NDSolve[\{f'[t] == h[t], 
    h'[t] == 12 f[t]*h[t] + 2 Coth[2 t]*h[t] - h[t]\^{}2/f[t], 
    f[-10\^{}(-5)] == 10\^{}(-5)/6.73, h[-10\^{}(-5)] == -1/6.73\}, \{f, 
    h\}, \{t, -100, -10\^{}(-5)\}];
    
    Plot[Evaluate[f[t] /. v], \{t, -40, -10\^{}(-5)\}, AxesLabel -$>$ \{r, f\}, 
 PlotRange -$>$ \{\{-30, -10\^{}(-5)\}, \{-0.1, 0.35\}\}, 
 PlotLabel -$>$  \{P\^{}2, 2/3  D\}]
 
 Plot[Evaluate[h[t] /. v], \{t, -8, -10\^{}(-5)\}, 
 AxesLabel -$>$ \{r, Subscript[f, r]\}, 
 PlotRange -$>$ \{\{-8, 0\}, \{-0.16, 0.01\}\}, PlotLabel -$>$ \{P\^{}2, 2/3 D\}]
 
 \item Figure \ref{tautmax}
 
 \begin{enumerate}
 \item Figure \ref{tautmax-a}
 
  Array[q, 200];
For[k = 1, k $<$ 201, k++, q[k] = 0];
For[i = 1, i $<$ 201, 
  i++, \{u = 
    NDSolve[\{f'[t] == h[t], 
      h'[t] == 12 f[t]*h[t] + 2 Coth[2 t]*h[t] - h[t]\^{}2/f[t], 
      f[-10\^{}(-5)] == 10\^{}(-5)/(0.5 + 0.1*i), 
      h[-10\^{}(-5)] == -1/(0.5 + 0.1*i)\}, \{f, h\}, \{t, -300, -0.01\}],
   For[j = 0, -300 + 0.01*j $<$ -0.01, j++, 
    q[i] = q[i] + Sqrt[-Evaluate[h[-300 + 0.01*j] /. u]]*0.01]\}
  ];
  
  ListLinePlot[\{Table[\{0.5 + 0.1*i, Re[Extract[q[i], 1]]\}, \{i, 1, 300\}],
   Table[\{0.5 + 0.1*i, Pi/2\}, \{i, 1, 300\}]\}, 
 PlotRange -$>$ \{\{0, 30\}, \{0.5, 1.8\}\}, PlotStyle -$>$ \{Blue, Red\}, 
 AxesLabel -$>$ \{$\backslash$[Tau], Subscript[t, max]\}]
 
 \item Figure \ref{tautmax-b}
 
 Array[p, 50]; Array[q, 50];
For[k = 1, k $<$ 51, k++, q[k] = 0];
For[i = 1, i $<$ 51, 
  i++, \{u = 
    NDSolve[\{f'[t] == h[t], 
      h'[t] == 12 f[t]*h[t] + 2 Coth[2 t]*h[t] - h[t]\^{}2/f[t], 
      f[-10\^{}(-5)] == 10\^{}(-5)/(100*i), h[-10\^{}(-5)] == -1/(100*i)\}, \{f, 
      h\}, \{t, -300, -0.01\}], p[i] = Evaluate[f[-300] /. u],
   For[j = 0, -300 + 0.005*j $<$ -0.01, j++, 
    q[i] = q[i] + Sqrt[-Evaluate[h[-300 + 0.005*j] /. u]]*0.005]\}
  ];
 
 ListLinePlot[\{Table[\{100*i, Re[Extract[q[i], 1]]\}, \{i, 1, 50\}], 
  Table[\{100*i, Pi/2\}, \{i, 1, 50\}]\}, 
 PlotRange -$>$ \{\{0, 5000\}, \{1.2, 1.8\}\}, PlotStyle -$>$ \{Blue, Red\}, 
 AxesLabel -$>$ \{$\backslash$[Tau], Subscript[t, max]\}]
 \end{enumerate}
 
 \item Figure \ref{convdata}
 
 n = 4; Array[p, \{n, 300\}] ; Array[q, \{n, 300\}];  Array[R, \{n, 
  300\}]; Array[c, \{n, 300\}]; 
For[i = 1, i $<$ 5, i++,
  \{s = 0, 
   u = NDSolve[\{f'[t] == h[t], 
      h'[t] == 12 f[t]*h[t] + 2 Coth[2 t]*h[t] - h[t]\^{}2/f[t], 
      f[-10\^{}(-5)] == 10\^{}(-5)/(10\^{}i), h[-10\^{}(-5)] == -1/(10\^{}i)\}, \{f, 
      h\}, \{t, -300, -0.0001\}],
   For[k = 0, -300 + k $<$ -0.01, k++,
    \{For[j = 0, j $<$ 1000, 
      j++, \{s = 
        s + Sqrt[-Evaluate[h[-300 + k + 0.001*j] /. u]]*0.001\}], 
     p[i, k + 1] = Re[s],
     q[i, k + 1] = Evaluate[f[-300 + k + 0.001*j] /. u],
     R[i, k + 1] = -Tanh[-300 + k + 0.001*j], 
     c[i, k + 1] = Evaluate[-h[-300 + k + 0.001*j] /. u]\}
    ]\}];
    
    ListPlot[
 Join[Table[
   Table[\{Extract[p[i, n], 1], Extract[q[i, n], 1]\}, \{n, 1, 300\}], \{i,
     1, 4\}], \{Table[\{Pi*i/40, (1 + Cos[Pi*i/20])/8\}, \{i, 1, 20\}]\}],
 PlotRange -$>$ \{\{0, 2\}, \{-0.1, 0.3\}\}, 
 PlotStyle -$>$ \{Blue, Green, Orange, Pink, Red\}, 
 PlotMarkers -$>$ Automatic, AxesLabel -$>$ \{t, f\}]
 
 ListLinePlot[
 Join[Table[
   Table[\{Extract[p[i, n], 1], R[i, n]\}, \{n, 1, 300\}], \{i, 1, 
    4\}], \{Table[\{Pi*i/20, 1\}, \{i, 1, 10\}]\}],
 PlotRange -$>$ \{\{0, 2\}, \{-0.1, 1.2\}\}, 
 PlotStyle -$>$ \{Blue, Green, Orange, Pink, Red\}, 
 PlotMarkers -$>$ Automatic, AxesLabel -$>$ \{t, a/b\}]
 
ListPlot[
 Join[Table[
   Table[\{Extract[p[i, n], 1], Extract[c[i, n], 1]\}, \{n, 1, 300\}], \{i,
     1, 4\}], \{Table[\{Pi*i/40, (1 - Cos[i*Pi/10])/32\}, \{i, 1, 20\}]\}], 
 PlotRange -$>$ \{\{0, 2\}, \{0, 0.13\}\},
  PlotStyle -$>$ \{Blue, Green, Orange, Pink, Red\}, 
 PlotMarkers -$>$ Automatic, AxesLabel -$>$ \{t, c\^{}2\}]
 
 \item Figure \ref{bubble}
 
 u = Table[\{u = 
     NDSolve[\{f'[t] == h[t], 
       h'[t] == 12 f[t]*h[t] + 2 Coth[2 t]*h[t] - h[t]\^{}2/f[t], 
       f[-10\^(-5)] == 10\^{}(-5)/(5000*i), 
       h[-10\^(-5)] == -1/(5000*i)\}, \{f, h\}, \{t, -300, -0.0001\}]\}, \{i, 
    1, 10\}];
    
    Plot[\{Table[Evaluate[f[t] /. Extract[u, i]]*5000*i, \{i, 1, 10\}], 
  Sinh[-t]\}, \{t, -12, 0\}, AxesLabel -$>$ \{r, f\}]
  
  Plot[\{Table[
   Evaluate[h[t] /. Extract[u, i]]*5000*i, \{i, 1, 10\}], -Cosh[
    t]\}, \{t, -12, -10\^{}(-5)\}, AxesLabel -$>$ \{r, Subscript[f, r]\}]
 
\end{enumerate}

\vspace*{5mm}

\noindent Chi Li\\
Department of Mathematics, SUNY at Stony Brook.\\
Email Address: chi.li@stonybrook.edu

\end{document}